\documentclass[11pt]{article}

\usepackage{amssymb}

\usepackage{graphicx}

\usepackage{amssymb}
\usepackage{latexsym}

\def\mathsf{\bf}

\def\d{\mathrm d}

\def\R{\mathbb{R}}
\def\Z{\mathbb{Z}}

\def\E{\mathrm E}

\def\P{\mathrm P}

\def\text{\mbox}

{\catcode `\@=11 \global\let\AddToReset=\@addtoreset}
\AddToReset{equation}{section}

\newtheorem{theorem}{Theorem}[section]

\newtheorem{proposition}[theorem]{Proposition}
\newtheorem{definition}[theorem]{Definition}

\newtheorem{remark}[theorem]{Remark}
\newtheorem{example}[theorem]{Example}

\def\mathsf{\bf}

\def\R{\mathbb{R}}
\def\Z{\mathbb{Z}}

\def\text{\mbox}

\def\d{\mathrm{d}}

\begin{document}

\title{A two-sample test for comparison of long
memory parameters}

\author{Fr\'ed\'eric Lavancier${}^{(1)}$, Anne Philippe${}^{(1)}$, Donatas Surgailis${}^{(2)}$ }

\date{\today \\ \small \textbf{1.} Laboratoire de Math\'ematiques Jean Leray,
Universit\'e de Nantes, France\\
\texttt{\{Frederic.Lavancier,Anne.Philippe\}@univ-nantes.fr} \\
\textbf{2.} Institute of Mathematics and Informatics, Vilnius, Lithuania
\\
 \texttt{sdonatas@ktl.mii.lt} \\
}

\maketitle

\begin{abstract}

We construct a two-sample test for comparison of long memory parameters based on ratios of two rescaled variance (V/S) statistics
studied in [Giraitis L., Leipus, R., Philippe, A., 2006. A test for
  stationarity versus trends and unit roots for a wide class of dependent
  errors. Econometric Theory 21, 989--1029]. The two samples have the same length and can be mutually independent or dependent.
  In the latter case, the test statistic is modified to make it asymptotically free of the long-run correlation coefficient between
  the samples. To diminish the sensitivity of the test on the choice of the bandwidth parameter, an adaptive formula
  for the bandwidth parameter is derived using the asymptotic expansion in [Abadir, K., Distaso, W., Giraitis, L., 2009.
  Two estimators of the long-run
variance: Beyond short memory. Journal of Econometrics 150, 56--70]. A simulation study shows that the above choice of bandwidth leads to a good size of our comparison test for most values of fractional and
ARMA parameters of the simulated series.

\end{abstract}
\begin{quote}

\end{quote}

\vskip2cm

\section{Introduction}

Long memory is one of the most widely discussed ``stylized facts'' of financial time series (see, e.g., Teyssi\`ere and Kirman~(2007)). In real data, long memory can be confused with short memory, unit roots, trends, structural changes, heavy tails and other features. Various tests for long memory have
been developed in the literature. See Lo~(1991), Kwiatkowski et al.~(1992), Robinson~(1994), Lobato and Robinson~(1998),
Giraitis et al.~(2003, 2006), Surgailis et al.~(2008). Most of these results pertain to the case of a single sample.

A natural extension of one-sample test about unknown long memory parameter $d$ is two-sample testing for comparison of respective memory parameters $d_1$ and $d_2$. In particular, such test can be useful for the memory
propagation (from durations to counts and realized volatility), question discussed in Deo et al.~(2009). Several studies
compare estimates of long memory parameter from different foreign exchange data and other sources (Cheung~(1993), Soofi et al.~(2006),
Casas and Gao~(2008)). Two-sample testing is also
related to the change-point problem
of the memory parameter discussed in Beran and Terrin~(1996), Horv\'ath~(2001).

The present paper constructs a test for testing the null hypothesis $d_1=d_2$ that long memory parameters
$d_i \in [0, 1/2) $ of two samples of length $n$, taken from respective stationary processes $X_i$, $i=1,2$, are equal, against the alternative
$d_1 \neq d_2$. The test statistic, $T_n$, is defined as a sum
\begin{eqnarray}
T_{n}&=&\frac{V_1/ S_{11,q}}{V_2/ S_{22,q}}\ + \ \frac{V_2/S_{22,q}}{V_1/ S_{11,q}}, \label{TT}
\end{eqnarray}
of two ratios of V/S, or rescaled variance, statistics $V_1/ S_{11,q}$ and $ V_2/S_{22,q}$ computed from
samples $(X_1(1), \dots, X_1(n)) $ and $(X_2(1), \dots, X_2(n))$. Here, $V_i$ is the empirical variance
of partial sums of $X_i$ and $S_{ii,q}$ is the Newey-West or HAC estimator of the long-run variance of $X_i$. The V/S statistic
was developed in Giraitis et al.~(2003, 2006) following the works of Lo~(1991) and Kwiatkowski et al.~(1992) on related R/S type
statistics. In particular, from Giraitis et al.~(2006) one easily derives the asymptotic null distribution $T$ of the statistic $T_n$
under the condition that the two samples are independent. It is also easy to show that for $d_1 \neq d_2$, one of the ratios
in (\ref{TT}) tends to infinity and the other one to zero, meaning that the test is consistent against the alternative $d_1 \ne d_2$.

However, independence of the two samples is too restrictive and may be unrealistic in financial data analysis since price movements of different assets are usually correlated and susceptible to common macroeconomic shocks. In order to eliminate the eventual dependence between samples, a modification
$\tilde T_n$ of (\ref{TT}) is proposed, which uses residual observations $(\tilde X_1(1), \dots, \tilde X_1(n)) $, obtained by
regressing partial sums of $X_1$ on partial sums of $X_2$. The modified statistic $\tilde T_n$ is shown to have
the same limit null distribution $T$ as if the two samples were independent.

It is well-known that a major problem in applications of the rescaled variance and related statistics
is the choice of the bandwidth parameter $q$. The present paper contributes to this problem by providing
an adaptive formula in (\ref{qopt}) for ``optimal'' $q$ which depends not only on the (common) memory parameter $d$ but also on the difference
between estimated short memory (AR) components of the spectrum of the sampled series. The derivation of the last result
uses the expansion of the HAC estimator in Abadir et al.~(2009). A simulation study confirms that
using this choice of bandwidth leads to a good size of our comparison tests for most values of fractional and
ARMA parameters of the simulated series.

Several interesting open problems were suggested by Referee and Associated Editor. The assumption of stationarity can be very restrictive
for applications.  We expect that our results can be extended to values of $d_i$ outside the interval $[0, 1/2)$.
Another possibility for future research is development
of similar procedures to test equal memory for more than two series.

The plan of the paper is as follows. Section \ref{test-construction} formulates the settings of the paper (Assumptions A$(d_1,d_2)$ and B$(d_1,d_2)$) and derives the limit of the test statistics $T_n$ and $\tilde T_n$ and the rejection region of the null hypothesis $d_1 = d_2$.
Assumption A$(d_1,d_2)$ guarantees the existence of long-run (co)variances and the consistency of the HAC estimators. Assumption B$(d_1,d_2)$
specifies the joint limit behavior of partial sums of $X_1$ and $X_2$ as given by bivariate fractional Brownian motion. The last process
is defined by means of stochastic integral representation as in Chung~(2002).  The test procedures are then presented in detail and a brief study focus on the asymptotic power of the tests $T_n$ and $\tilde T_n$. Section \ref{bivariate-linear} verifies Assumptions A$(d_1,d_2)$ and B$(d_1,d_2)$ for bivariate moving average $(X_1,X_2)$.
Section \ref{simulation-study} assesses the performance of the tests $T_n$ and $\tilde T_n$, by simulating bivariate FARIMA samples with various fractional and autoregressive/moving average parameters. Conclusions are given in Section \ref{conclusion}. The Appendix contains auxiliary results and derivations.

\smallskip

\noindent {\bf Notation.} Below, $\rightarrow_p, \rightarrow_{\rm law}, \rightarrow_{D[0,1]}$ and $\rightarrow_{\rm fdd} \ (=_{\rm fdd})$ stand for the convergence in probability, the weak convergence of random variables, the weak convergence of random elements in the Skorohod space
$D[0,1]$, and the weak convergence (equality) of finite dimensional distributions, respectively. Relation `$\sim$' means that the ratio
of both sides tends to 1.

\section{Construction of tests and its properties}\label{test-construction}

\subsection{Assumptions and main results}

Let $((X_1(t), X_2(t)), \ t \in {\Z})$ be a bivariate covariance stationary process, viz.,
a sequence of random vectors $(X_1(t), X_2(t)) \in {\R}^2 $ such that ${\E} X_i(t) = \mu_i $ and
$$
{\rm cov}(X_i(t), X_j(t+h)) = \gamma_{ij}(h)
$$
do not depend on $t\in {\Z}$ for any $h \in \Z, i,j=1,2$. Introduce
the popular Bartlett-kernel estimator of the long-run (co)variance:
\begin{eqnarray}
S_{ij,q}&=&\sum_{h=-q}^q \left(1-\frac{|h|}{q+1}\right)\hat
\gamma_{ij}(h), \label{Sdef}
\end{eqnarray}
where
\begin{eqnarray} \label{hatgamma}
\hat \gamma_{ij}(h)&=&n^{-1}\cases{\displaystyle\sum_{t=1}^{n-h} (X_i(t) - \bar
X_i)(X_j(t+h)- \bar X_j),&$h\ge 0$, \cr \displaystyle\sum_{t=1-h}^{n} (X_i(t) - \bar
X_i)(X_j(t+h)- \bar X_j),&$h\le 0$,\cr}
\end{eqnarray}
$\bar X_i = n^{-1} \sum_{t=1}^n X_i(t)$.
The estimator in (\ref{Sdef}) is also called the heteroskedasticity and autocorrelation consistent (HAC).
Also define
\begin{eqnarray}
S^{\circ}_{ij,q}&=&\sum_{h=-q}^q \left(1-\frac{|h|}{q+1}\right)\hat
\gamma^\circ_{ij}(h), \label{Scirc}
\end{eqnarray}
where
\begin{eqnarray} \label{circgamma}
\hat \gamma^\circ_{ij}(h)&=&n^{-1}\cases{\displaystyle\sum_{t=1}^{n-h} (X_i(t) -
\mu_i)(X_j(t+h)- \mu_j),&$h\ge 0$, \cr \displaystyle\sum_{t=1-h}^{n} (X_i(t) -
\mu_i)(X_j(t+h)- \mu_j),&$h\le 0$.\cr}
\end{eqnarray}
Note $\hat \gamma_{ij}(h) = \hat \gamma_{ji}(-h), \ \hat \gamma^\circ_{ij}(h) = \hat \gamma^\circ_{ji}(-h), \
S_{12,q} = S_{21,q}, \ S_{12,q}^\circ = S_{21,q}^\circ.$

\bigskip

\noindent \underline{\sc Assumption A$(d_1,d_2)$}\ \ There exist $d_i \in [0,1/2), i =1,2 $ such that for any $i,j=1,2$ the following limits exist
\begin{eqnarray}
c_{ij}&=& \lim_{n\to \infty} \frac{1}{n^{1+d_i +d_j}} {\E} \left(\sum_{t=1}^n (X_i(t)- \mu_i)\right)\left(\sum_{s=1}^n (X_j(s)- \mu_j)\right)\nonumber \\
&=&\lim_{n \to \infty} \frac{1}{n^{1+d_i+d_j}} \sum_{t,s=1}^n \gamma_{ij}(t-s).\label{Clim1}
\end{eqnarray}
Moreover,
\begin{eqnarray}
\frac{\displaystyle\sum_{k,l=1}^q \hat \gamma_{ij}(k-l)}{\displaystyle
\sum_{k,l=1}^q \gamma_{ij}(k-l)}
&\rightarrow_p&1 \label{Clim2}
\end{eqnarray}
as $q \to \infty, n \to \infty, n/q \to \infty$.

\begin{remark} {\rm The asymptotic constant $c_{ij}$ is called the long-run covariance of $X_i$ and $X_j$. Condition (\ref{Clim2}) is
similar to Giraitis et al.~(2006, Assumption A.2). It guarantees the consistency of the HAC estimator (see below).
}
\end{remark}

\begin{proposition} \label{propS} Let Assumption A$(d_1,d_2)$ hold. Then, as $q \to \infty, n \to \infty, n/q \to \infty$,
\begin{eqnarray}
q^{-d_i-d_j} S_{ij,q}&\rightarrow_p&c_{ij}, \quad q^{-d_i-d_j} S^\circ_{ij,q}\ \rightarrow_p\ c_{ij} \quad (i,j=1,2). \label{Slimit}
\end{eqnarray}
Moreover,
\begin{equation} \label{SS}
\frac{1}{q}\big(S_{ij,q} - S^\circ_{ij,q}\big) \ = \ - (\bar X_i - \mu_i) (\bar X_j- \mu_j)\Big(1 + o_p(1)\Big).
\end{equation}

\end{proposition}

Assumption B$(d_1,d_2)$, below, specifies the joint limit of partial sums of $X_1$ and $X_2$. It is similar to Giraitis et al.~(2006,
Assumption A.1). The limit process (bivariate fractional Brownian motion)  is defined
 through a stochastic integral representation (\ref{Bi}) similarly as in Chung~(2002, (6)). Equivalently, this process can be defined
 through the covariance function defined in   (\ref{covB12}).

\begin{definition} \label{bifBm} A  nonanticipative bivariate fractional Brownian
motion (bi-fBm) with memory
parameters $d_i \in (-1/2, 1/2), \ i=1,2,$ is a Gaussian process $B = \big( (B_{1}(s), B_{2}(s)), s \in \R) $
admitting the following representation for $i=1,2$
\begin{equation}
B_{i} (t)\ = \ \cases{c(d_i) \int_{\R} \left((t -x)_+^{d_i} -
(-x)_+^{d_i }\right) W_i(\d x), &if $d_i \ne 0$, \cr
W_i(0,t], &if $d_i = 0$, \cr}
\label{Bi}
\end{equation}
$W = ((W_1({\d}x), W_2({\d}x)),\ x \in \R )$
is a $2-$dimensional Gaussian independently scattered white noise
with real components, zero mean and covariance matrix
\begin{equation}\label{Wij}
{\E} W_i({\d}x)W_j(\d x) = {\d}x \cases{1, &$i=j$, \cr
\rho_W, &$i\ne j$\cr},
\end{equation}
for some $\rho_W \in [-1,1]$,
and the constants $c(d_i) $ are determined by condition \ $\E B^2_i(1) =1 $ so that
$$
c^2(d_i) = \Big(\int \big((1 -x)_+^{d_i} -
(-x)_+^{d_i }\big)^2 \d x\Big)^{-1} = \frac{\cos (d_i \pi)}{{\rm B}(d_i+1,d_i+1)},
$$
where ${\rm B}(p,q)$ is the beta function and $x_+ = \max(x, 0)$.
\end{definition}

\begin{remark} \label{rho} {\rm The nonanticipative bi-fBm is a particular case of general bi-fBm having the stochastic representation
$$
X(t)\ = \ \int_{\R} \left\{ \left( (t-x)_+^{D}- (-x)_+^{D}\right) A_+ +
\left((t-x)_-^{D} - (-x)_-^{D}\right)A_- \right\} \widetilde W({\rm d}x),
$$
where $D= {\rm diag}(d_1,d_2), \,  x_- = \max(-x,0)$, $A_+,A_-$ are real $2\times 2$ matrices and $\widetilde W({\rm d}x)= (\widetilde W_1({\rm d}x),\widetilde W_2({\rm d}x)), \, x \in \R$, is an independently scattered Gaussian white noise with zero mean, unit variance and independent components;
see Didier and Pipiras (2010),  also  Lavancier et al.~(2009, (1.6)).
The representation in (\ref{Bi}) corresponds to the matrices
 $A_+A_+^* = \pmatrix{ c(d_1)^2 &c(d_1)c(d_2)\rho \cr
 c(d_1)c(d_2)\rho &c(d_2)^2 \cr}, $  $\
A_-=0$.}
\end{remark}

\begin{remark} \label{cov}
{\rm In the sequel by bi-fBm we mean the nonanticipative process in (\ref{Bi}).
A bi-fBm has stationary increments and the self-similarity property:
$$(\lambda^{-d_1 - .5} B_1(\lambda t), \lambda^{-d_2-.5} B_2(\lambda t)) \ =_{\rm fdd} \
(B_1(t), B_2(t))$$ for any $\lambda>0$. Lavancier et al.~(2009)  showed that these properties essentially determine
the covariance function in (\ref{covB12})-(\ref{covB12L}), up to a choice of constants $g_{ij}$ and $ g_i $ defined in (\ref{eq:2}) and (\ref{eq:3}) (see Section \ref{sec:covar-funct-bivar}).

Note that each component $B_i$ is a univariate fractional Brownian motion with variance $\E B_i^2(t) = |t|^{2d_i+1}, \ i=1,2$.
Also
note that $(B_1(1), B_2(1)) $ has a bivariate Gaussian distribution with zero means, unit variances and the correlation
coefficient $\rho = \E B_1(1)B_2(1) = \rho_W \kappa(d_1,d_2)$, where $\kappa(d_1,d_2)$ depends only on $d_1,d_2$. In the case
$d_1=d_2$, we have  that $\rho = \rho_W $ and the process $(\tilde B_1(t) = B_1(t) - \rho B_2(t), t \in \R) $ is a fractional Brownian
motion with variance $\E \tilde B^2_1(t) = (1- \rho^2)|t|^{2d_1+1} $. Moreover, the processes $\tilde B_1$ and $B_2$ are {\it independent}.
Indeed, from (\ref{Bi}) it is immediate that $\tilde B_1$ has a similar stochastic representation with $W_1({\d}x)$ replaced
by $\tilde W_1({\d}x) = W_1({\d}x) - \rho W_2({\d}x)$, with $\tilde W_1 $ independent of $W_2$.

}
\end{remark}

\bigskip

\noindent \underline{\sc Assumption B$(d_1,d_2)$} \ \ Assumption A$(d_1,d_2)$ is satisfied and, moreover,
\begin{eqnarray}
&&\left(n^{-d_1-(1/2)}\sum\nolimits_{t=1}^{[n\tau]} (X_1(t)- {\E}X_1(t)),
n^{-d_2-(1/2)}\sum\nolimits_{t=1}^{[n\tau]} (X_2(t)- {\E} X_2(t)) \right) \nonumber \\
&&\rightarrow_{\rm fdd} \
\left(\sqrt{c_{11}} B_{1}(\tau), \sqrt{c_{22}} B_{2}(\tau)\right), \label{biconv}
\end{eqnarray}
where $c_{ij}$ are the same as in (\ref{Clim1}), $c_{ii} >0, \ i=1,2$ and $(B_{1}, B_{2})$ is a bi-fBm
with memory parameters $d_1, d_2$ and the correlation coefficient $\rho = {\rm corr}(B_1(1), B_2(1)) = c_{12}/\sqrt{c_{11} c_{22}} \in [-1,1]$.

\bigskip

Define the empirical variance of partial sums of $X_i$:
\begin{eqnarray}
V_i&=&
n^{-2}\sum^n_{k=1}\left(\sum^k_{t=1}(X_i(t)-\bar X_i)\right)^2- n^{-3}\left(\sum^n_{k=1}\sum^k_{t=1}(X_i(t)-\bar X_i)\right)^2.
\label{V}
\end{eqnarray}
The following proposition obtains the limit distribution of the statistic $T_n$ defined in (\ref{TT}).

\begin{proposition} \label{Test} Let Assumptions A$(d_1,d_2)$ and B$(d_1,d_2)$ be satisfied
with some $d_1, d_2 \in [0,1/2)$ and $\rho \in [-1,1]$, and let $n, q, n/q \to \infty$.

\smallskip

\noindent (i) If $d_1 = d_2 = d $ then
\begin{eqnarray}
T_n&\rightarrow_{\rm law}&T \ = \ \frac{U_1}{U_2} + \frac{U_2}{U_1}, \label{Tlimit1}
\end{eqnarray}
where
\begin{eqnarray}
U_i&=&
\int_0^1 (B^0_i(\tau))^2 {\d}\tau
- \left(\int_0^1 B^0_i(\tau) {\d}\tau\right)^2, \quad i=1,2,\label{Uidef}
\end{eqnarray}
where
$(B^0_i(\tau) = B_i(\tau) - \tau B_i(1), \tau \in [0,1]), \ i=1,2 $ are fractional Brownian bridges obtained
from bivariate fBm $((B_1(\tau), B_2(\tau)), \tau \in \R)$ with the same memory parameters $d_1 = d_2 = d$ and
correlation coefficient $\rho = \rho_W $ (see Definition \ref{bifBm}).

\smallskip

\noindent (ii) If $d_1 \neq d_2$ then
\begin{equation}
T_n \ \rightarrow_p \ \infty.
\end{equation}

\end{proposition}

Let $\tilde V_1, \tilde S_{11,q} $ be the statistics in (\ref{V}), (\ref{Sdef}),
respectively, where $X_1(t), t=1,\cdots, n$ is replaced by
\begin{equation}
\tilde X_1(t) \ = \ X_1(t) - (S_{12,q}/S_{22,q}) X_2(t), \qquad t=1,\dots, n. \label{TildeX}
\end{equation}
In particular, note
\begin{eqnarray}
\tilde S_{11,q}&=&S_{11,q} - \frac{S^2_{12,q}}{S_{22,q}}. \label{tildeS}
\end{eqnarray}
Define
\begin{eqnarray} \label{tildeT}
\tilde T_n&=&
\frac{\tilde V_1/ \tilde S_{11,q}}{V_2/ S_{22,q}}\ + \ \frac{V_2/S_{22,q}}{\tilde V_1/ \tilde S_{11,q}}.
\label{tildeTn}
\end{eqnarray}
Note, $\tilde T_n$ is obtained by replacing $V_1, S_{11,q}$ in the definition of $T_n$ in (\ref{TT}) by the
corresponding quantities $\tilde V_1, \tilde S_{11,q}$ as defined above.

In the following proposition, we prove that under the null hypothesis, the limit distribution of $\tilde T_n$ is free of $\rho$,  contrary to $T$ in  (\ref{Uidef}). Note that the limit of $\tilde T_n$ coincides with (\ref{Uidef}) when $\rho=0$.   This occurs for example when the statistics $T_n$  is calculated from two independent processes  $X_1,\; X_2 $.

\begin{proposition} \label{TTest} Let Assumptions A$(d_1,d_2)$ and B$(d_1,d_2)$ be satisfied
with some $d_1, d_2 \in [0,1/2)$ and $\rho \in (-1,1)$, and let $n, q, n/q \to \infty$.

\smallskip

\noindent (i) If $d_1 = d_2 = d $ then
\begin{eqnarray}
\tilde T_n&\rightarrow_{\rm law}&\tilde T \ = \ \frac{\hat U_1}{U_2} + \frac{U_2}{\hat U_1},
\label{Tlimit}
\end{eqnarray}
where $\hat U_1, U_2$ are independent and have the same distribution in (\ref{Uidef}).

\smallskip

\noindent (ii) If $d_1 > d_2$ then
\begin{equation} \label{Tinfty}
\tilde T_n \ \rightarrow_p \ \infty.
\end{equation}

\smallskip

\noindent (iii) If $d_1 < d_2 $ then
\begin{equation} \label{Trho}
\tilde T_n \ \rightarrow_p \ \frac{\rho^2}{1-\rho^2} + \frac{1-\rho^2}{\rho^2}.
\end{equation}

\end{proposition}

\begin{remark} {\rm The ratio $\hat \beta = S_{12,q}/S_{22,q} $ in (\ref{TildeX})
minimizes the sum of squares:
$$
\sum_{i=1-q}^n \Big(\sum_{t=i\vee 1}^{(i+q)\wedge n}(X_1(t)- \bar X_1)-
\beta \sum_{t=i\vee 1}^{(i+q)\wedge n}(X_2(t)- \bar X_2)\Big)^2.
$$
Therefore, $\hat \rho = \hat \beta \sqrt{S_{22,q}/S_{11,q}} $ can be considered as the least squares estimate of the long-run correlation coefficient
$\rho $  between partial sums
of the two samples.

}
\end{remark}

\subsection{Testing procedures}
\label{sec:testing-procedures}

Let $t_{\alpha} (d) $ denote the upper $\alpha-$quantile of the
r.v. $\tilde T$  defined in (\ref{Tlimit}) (or $T$ in (\ref{Uidef})  when $\rho=0$), viz.
\begin{equation} \label{critic}
\alpha = \P (\tilde T> t_\alpha(d)), \qquad \alpha \in (0,1).
\end{equation}
Let
\begin{equation}
\hat d = (\hat d_1 + \hat d_2)/2, \label{eq:1}
\end{equation}
where $\hat d_i $ is an estimator of $d_i$ satisfying
\begin{equation}
\hat d_i - d_i = o_p(1/\log n) \qquad (i=1,2).
\end{equation}
Similarly to Giraitis et al.~(2006, Lemma 2.1), it can be proved that
the quantile function $t_{\alpha} (d)$ is continuous in $d \in
[0,1/2)$ for any $\alpha \in (0,1)$. Therefore, the estimated quantile
$t_{\alpha} (\hat d) \rightarrow_p t_{\alpha} (d) $ as $n \to \infty $
and the asymptotic level of the tests associated to the critical regions in (\ref{Trule})-(\ref{TTrule}) is preserved by replacing $t_{\alpha}
(d)$ by $t_{\alpha} (\hat d) $.

\vskip .5cm

\noindent{\it Testing the equality of the memory parameters in the case of independent samples.}\/ We wish to test the
null hypothesis  $d_1 = d_2 $ against the alternative $d_1 \ne 
d_2$ under the assumption that $X_1$ and $X_2$ are independent.  The decision rule at $\alpha-$level of significance is the
following:  we reject the null hypothesis
 when
\begin{equation} \label{Trule}
T_n > t_\alpha(\hat d).
\end{equation}
The consistency of this test is ensured by Proposition \ref{Test} (ii).

\vskip .5cm

\noindent{\it Testing the equality of the memory parameters in the case of possibly dependent samples.}\/ We wish to test the
null hypothesis  $d_1 = d_2 $ against the alternative $d_1 > 
d_2$ in the general case when  $X_1$ and $X_2$ are possibly dependent.  The decision rule at $\alpha-$level of significance is the
following:  we reject the null hypothesis
 when
\begin{equation} \label{TTrule}
\tilde T_n > t_\alpha(\hat d).
\end{equation}
The consistency of this test is ensured by Proposition \ref{TTest} (ii).

\begin{remark} {\rm 
  For testing $d_1 = d_2$ against $d_1 <d_2$, the samples $X_1$ and
  $X_2$ should be exchanged in the statistic (\ref{tildeT}). }
\end{remark}

\begin{remark} {\rm As noted by the referee, an undesirable feature of the testing procedure in (\ref{TTrule}) is that a rejection might occur not only when $d_1>d_2$ but also when $d_1<d_2$ (see Proposition \ref{TTest} (iii)). To alleviate this feature, in (\ref{TTrule}) one can use the statistic
$$
\tilde T_n^+ \ = \  \frac{\tilde V_1/ \tilde S_{11,q}}{V_2/ S_{22,q}}
$$
instead of $\tilde T_n$.  Note $\tilde T_n = \tilde T_n^+ + (\tilde T^+_n)^{-1}$. Under the assumptions of Proposition \ref{TTest}, the limit
distribution  of $\tilde T_n^+$ can easily be  obtained from the
proof of this proposition. In particular, 
if $d_1=d_2=d$, then
$$
\tilde T_n^+ \ \rightarrow_{\rm law}\  \tilde T^+ = \frac{\hat U_1}{U_2},
$$
where $\hat U_1, U_2$ are independent and have the same distribution
in (\ref{Uidef}). From the proof of Proposition \ref{TTest}, it also follows easily that 
$\tilde T_n^+ \to_p \infty \ (d_1 > d_2)$ and $\tilde T_n^+ \to_p \rho^2/(1-\rho^2)\  (d_1<d_2)$, i.e., $\tilde T_n^+$ does not explode to infinity 
for $d_1<d_2$ and $|\rho|$ not very close to 1.}
\end{remark}



\subsection{Asymptotic behavior of the power function}
We discuss in this section the asymptotic power of the testing procedures in (\ref{Trule}) and (\ref{TTrule}). 

For testing $d_1=d_2$ against $d_1>d_2$,  we have proved that the rejection
probability of the null hypothesis  tends to 1, i.e., that
\begin{equation} \label{pow1}
\P (\tilde T_n > t_\alpha(\hat d)) \to 1,
\end{equation}
for any $\alpha \in (0,1)$.

 Section~\ref{simulation-study} provides finite sample rejection frequencies
of the null hypothesis for some choices of parameters $d_1, d_2$ and some bivariate FARIMA  models. A natural question in this context
is to estimate the convergence   rate in (\ref{pow1}), or the decay rate of the probability
$\P (\tilde T_n \le a) $, as a function of $a, n, d_1, d_2, \rho $ and (possibly)  other quantities of the model assumptions.

From the proof of Proposition \ref{TTest} (ii), see (\ref{d12}) below, we have that for $d_1 > d_2$, the normalized
statistic  $\tilde T_n $ has a limit distribution,
\begin{equation} \label{pow2}
(q/n)^{2(d_1-d_2)} \tilde T_n  \ = \  (q/n)^{2(d_1-d_2)} \frac{\tilde V_1/ \tilde S_{11,q}}{V_2/ S_{22,q}} + o_p(1) \
\rightarrow_{\rm law} \  \frac{U_1}{(1-\rho^2)U_2} = \Lambda,
\end{equation}
say, where $U_i, i=1,2$ are defined in (\ref{Uidef}). Therefore, we can expect that
the probability $\P (\tilde T_n \le a)$ decays as the probability $\P (\Lambda \le a (q/n)^{2(d_1-d_2)}) $, when $n \to \infty$.
The decay rate of $\P( \Lambda \le x), \ x \to 0$ is unknown, even for independent $U_1, U_2$, but in principle can be estimated
from Monte-Carlo experiments. It is also plausible that r.v. $\Lambda $ has a bounded probability density near $x=0$
and so the above discussion suggests a decay rate  $\P (\tilde T_n \le a) = O( (q/n)^{2(d_1-d_2)})$.

However, the above argument is heuristic; in particularly, the replacement of $\tilde T_n$ by $(q/n)^{2(d_2-d_1)} \Lambda $ is
not rigorously justified.  It is clear that in order to correctly assess the probability $\P (\tilde T_n \le a)$, it is necessary
to control from above the probabilities of {\it high} values of  $\tilde V_1 $ and $S_{22,q} $ in the denominator of the ratio
$\frac{\tilde V_1/ \tilde S_{11,q}}{V_2/ S_{22,q}} = \frac{\tilde V_1 S_{22,q}}{V_2 \tilde S_{11,q}}$
 and the probabilities
of {\it small} values of $V_2,  \tilde S_{11,q} $ in the numerator of the last ratio. While the former probabilities can
be controlled by the Markov inequality, direct estimation of the latter probabilities is  difficult and is replaced by an
assumption on the distribution functions of
these r.v.'s, see  (\ref{A0}) below.

\begin{proposition} \label{power}
Let Assumption A$(d_1,d_2)$ be satisfied with $0\le d_2 < d_1 < 1/2$. Moreover, assume that the distribution functions of the normalized r.v.'s
$n^{-2d_1} V_1$ and $q^{-2d_2} S_{22,q}$ satisfy the following condition: for any $M >0, a_0>0$  there exists a constant $K$
such that the inequalities
\begin{equation} \label{A0}
\P(n^{-2d_1} V_1 \le a) \le Ka, \qquad \P(q^{-2d}S_{22,q} \le a) \le Ka
\end{equation}
are satisfied for any $n >M, q> M, n/q > M$ and any $0 < a \le a_0$. Then there exists a constant $K_1$, independent of $n, q, a $, and such
that
\begin{eqnarray} \label{A00}
\P (T_n \le a)
&\le&K_1 a^{1/4} \left(\frac{q}{n}\right)^{(d_1-d_2)/2}, \qquad
\P (\tilde T_n \le a)
\  \le \ K_1 a^{1/4} \left(\frac{q}{n}\right)^{(d_1-d_2)/2}
\end{eqnarray}
hold for all $n, q, n/q$ sufficiently large and any $a \ge 0$ from a compact set.
\end{proposition}

The proof of the above proposition is given in Section \ref{sec:proofs-propositions}.
Note that condition (\ref{A0}) is implied by the existence of bounded probability densities of  $n^{-2d_1} V_1$ and $q^{-2d_2} S_{22,q}$.
Also note that the assumptions of Proposition \ref{power} refer to {\it marginal} distributions of $V_1$ and $S_{22,q}$ only, and do not
impose a restriction on the joint distribution of the four statistics in the definition of $T_n $ and $\tilde T_n$.

\section{Application to bivariate linear processes with long memory.}\label{bivariate-linear}

In this section we specify Assumptions A$(d_1,d_2)$ and B$(d_1,d_2)$ to a class of bivariate linear models $(X_1(t),X_2(t)), \ t\in {\Z}$
as given by
\begin{eqnarray}
X_i(t)&=&\sum_{k=0}^\infty \psi_{i1}(k) \xi_1(t-k) + \sum_{k=0}^\infty \psi_{i2}(k) \xi_2(t-k), \qquad i=1,2, \label{Lin}
\end{eqnarray}
where $\psi_{ij}(k)$ are real coefficients with
$\sum_{k=0}^\infty \psi^2_{ij}(k) < \infty $ and
$(\xi_1(t), \xi_2(t)),$ $ t\in {\Z}$ is a bivariate (weak) white noise with nondegenerate covariance matrix $(\rho_{\xi,ij})_{i,j=1,2}$. In other words, $(\xi_1(t), \xi_2(t)),$ $ t\in {\Z}$ is a
sequence of random vectors
with zero mean ${\E} \xi_1(t) = {\E} \xi_2 (t) = 0$ and covariances
\begin{equation}
{\E} \xi_i(t) \xi_j(s) \ = \ \cases{\rho_{\xi,ij}, &$t=s$, \cr
0, &$t\ne s$. \cr} \label{cov}
\end{equation}
Without loss of generality, below we shall assume $\rho_{\xi,11} = \rho_{\xi, 22} = 1, \ \rho_{\xi, 12} = \rho_{\xi, 21} = \rho_\xi \in (-1,1)$.
\bigskip

\noindent \underline{{\sc Assumption} $\widetilde {\rm A}(d_{ij})$}\ $(X_1(t),X_2(t))$ is a bivariate linear covariance stationary process
as in (\ref{Lin}) with coefficients $\psi_{ij}(k)$ satisfying the following conditions:
\begin{itemize}
\item If $ d_{ij} \in (0,1/2) $
$$\psi_{ij}(k)\ = \  \left(\alpha_{ij} + o(1)\right) |k|^{d_{ij} - 1} \quad (k\to \infty) $$
where $\alpha_{ij} \ne 0 $ are some numbers, $i,j=1,2$.
\item If $ d_{ij}=0$
$$\sum_{k=0}^\infty |\psi_{ij}(k)| < \infty, \qquad \alpha_{ij} = \sum_{k=0}^\infty \psi_{ij}(k).$$
\end{itemize}

\bigskip

\noindent \underline{{\sc Assumption} $\widetilde {\rm B}(d_{ij})$}\ \ Assumption $\widetilde {\rm A}(d_{ij})$
is satisfied and, moreover, $(\xi_1(t),\xi_2(t)),$ $ t\in {\Z}$ is a sequence of i.i.d. random vectors.

\begin{proposition} \label{prop-bilinear}

(i) Let $(X_1(t),X_2(t))$ satisfy Assumption $\widetilde A(d_{ij})$. Then the limits $c_{ij}$ in Assumption A$(d_1,d_2)$,
(\ref{Clim1}) exist, with
\begin{equation}
d_i = \max \{d_{i1}, d_{i2}\} \in [0,1/2) \qquad (i=1,2). \label{di}
\end{equation}

\smallskip

\noindent (ii) Let $(X_1(t),X_2(t))$ satisfy Assumption $\widetilde B(d_{ij})$ and there exists $\delta> 0$ such that $ {\E} |\xi_i(t)|^{2+\delta} < \infty
\ (i=1,2)$. Then $(X_1(t),X_2(t))$ satisfies Assumptions $A(d_1,d_2)$ and $B(d_1,d_2)$, with $d_i$ as defined in (\ref{di}). Moreover,
the finite-dimensional convergence in Assumption $B(d_1,d_2)$, (\ref{biconv}) extends to the functional
convergence in the Skorohod space $D[0,1]$.

\end{proposition}

\begin{remark}

 {\rm Proposition \ref{prop-bilinear} (ii) complements Chung~(2002, Th.1), who discussed convergence
of partial sums of $K-$variate linear processes to $K-$variate fractional Brownian motion under slightly different
assumptions. Proposition \ref{prop-bilinear} (i) and Proposition \ref{propS} (\ref{Slimit}) also complement the result in Abadir et al.~(2009)
about consistency of the HAC estimator
for linear processes, by relaxing the 4th moment condition on the noise in the case $d_i \in [1/4,1/2)$.

}

\end{remark}

Let us consider some parametric examples of bivariate linear processes.
Hereafter, we denote by $L$ the backward shift  i.e.
$L X(t) = X(t-1)$.

\begin{example} \label{ex1}{\rm
Let $a_{ij} \in {\R} \ (i=1,2)$ be some constants, and let
\begin{eqnarray} \label{modeldep}
X_i(t)&=&(1-L)^{-d_i} (a_{i1} \xi_1(t) + a_{i2} \xi_2(t)) \qquad (i=1,2)
\end{eqnarray}
be FARIMA$(0,d_i,0)$ processes, with possibly different
parameters $d_i \in (0,1/2)$.
This process satisfies Assumption $\widetilde {\rm A}(d_{ij})$ with $d_{ij} = d_i$ and
$\alpha_{ij} = \Gamma(d_i)^{-1}a_{ij} \ (i,j=1,2)$.

 If $(\xi_i(t), \xi_2(t))$ form a sequence
of i.i.d. vectors as in Assumption $\widetilde {\rm B}(d_{ij})$, partial sums of $(X_1(t), X_2(t))$ converge
to a bivariate fBm $(B_{1,d_1}, B_{2,d_2})$.

The limiting fBm has independent components if and only if the noise has uncorrelated components, i.e., if \
${\E} (a_{11}\xi_1(t) + a_{12}\xi_2(t))(a_{21} \xi_1(t) + a_{22}
\xi_2(t))= 0$. For $d_1=d_2$, the last condition is equivalent to the uncorrelatedness of the components of the process:
${\E} X_1(t_1) X_2(t_2) = 0 \ (t_1,t_2 \in {\Z})$.
}
\end{example}


\begin{example} {\rm
Consider the following system of linear equations:
\begin{eqnarray*}
(1-L)^{d'_{11}} X_1(t) + \beta (1-L)^{d'_{12}}X_2(t)&=&\xi_1(t), \\
(1-L)^{d'_{22}}X_2(t)&=&\xi_2(t),
\end{eqnarray*}
where $d'_{ij} \in [0,1/2), \beta \in {\R}$ are parameters, $d'_{22}+ d'_{12} - d'_{12} < 1/2$ and where $(\xi_1(t),\xi_2(t)), t\in \Z $ are as in
(\ref{Lin}). A covariance stationary
solution of the above equation is given by
\begin{eqnarray*}
X_2(t)&=&(1-L)^{-d'_{22}}\xi_2(t), \\
X_1(t)&=&(1-L)^{-d'_{11}}\xi_1(t)- \beta (1-L)^{d'_{12}-d'_{11}- d'_{22}} \xi_2(t). \label{ex2}
\end{eqnarray*}
Then $(X_1(t), X_2(t))$ satisfies Assumption $\widetilde {\rm A}(d_{ij})$ with $d_{11} = d'_{11}, d_{12} =
d'_{22} + d'_{11}- d'_{12}, d_{21} = 0 $ and $d_{22} = d'_{22}$.

Assume $(\xi_1(t),\xi_2(t)) $ is a sequence of i.i.d. vectors satisfying the conditions in Proposition \ref{prop-bilinear} (ii) and
let $\beta \ne 0$. There are three cases $d'_{22} > d'_{12} $, $d'_{22} < d'_{12} $ and $d'_{22} = d'_{12}$ leading
to  $d_{11} < d_{12},\  d_{11} > d_{12}$ and $d_{11} = d_{12}$, respectively. In each of these cases we can determine the memory parameter
$d_i$ of $(X_i(t)), i=1,2$ and the limiting bi-fBm in Assumption B$(d_1,d_2)$, together with the correlation coefficient $\rho$.
}
\end{example}

\section{A simulation study}\label{simulation-study}

In this section we assess the finite-sample performance of our
procedures to test $d_1=d_2$ versus $d_1>d_2$
and provide a practical recommendation for the choice of the bandwidth parameter $q$.

 The memory parameters
$d_1$ and $d_2$ are estimated with the help of the adaptive version of the FEXP estimator (see Iouditsky {\it et al.} (2001)),
which in practice turns out to be  less sensitive to the short memory part of the long memory process as compared to other estimators.
The bandwidth parameter is chosen
according to the adaptive formula (\ref{qopt}) derived in  Appendix (see Section \ref{sec:deriv-adapt-bandw}). The optimisation of the bandwidth is realised under the null hypothesis in order to ensure a good size to the test procedure.  The choice of $q$ depends on $\hat d=\frac 12 (\hat d_1+\hat d_2)$
as in the expansion obtained by  Abadir et al.~(2009, (2.14)), but also takes into account the short memory spectrum, in a form of certain coefficient.

The simulated samples are independent or dependent Gaussian
FARIMA processes with different fractional and autoregressive/moving average parameters.
The 5\% quantile function
in (\ref{critic}) was approximated from extensive Monte-Carlo experiments by
\begin{equation}\label{pol-quantile}
t_{5\%} (d) \approx 3.7d^2 + 8.6d + 5.2.
\end{equation}

\noindent \textit{Independent samples.}\/
Table 1-3  concerns the case of independent samples, the  test
procedure is based on  $T_n$.
Tables 1 and 2 report the percentages of rejection of the null hypothesis $d_1=d_2$ of the test $T_n > t_{5\%}(\hat d)$ from 1000 replications
of {\it independent} FARIMA$(1,d,0)$ samples of size $n\in\{1024, \ 4096\}$, for five values of $d_i \in \{0,\ .1, \ .2, \ .3, \ .4\}$
and three values $a_i \in \{0, \ .4, \ .8 \}$ of the autoregressive parameter. Recall from Proposition \ref{Test} that for independent
samples $T_n$ and $\tilde T_n$ have the same limiting distribution.
We can see from both tables that the $T_n$ test has fairly good size for most values of $d_i$ and $a_i$. Table 3 provides the mean values of $\hat q$. Since these values are rather scattered across the table and
the rejection frequency is very sensitive to the choice of $q$, the general impression from Tables 1-3
is that the adaptive
choice of $q$ in (\ref{qopt}) was necessary. We also note that the
power of the test decreases with increase of $|a_1-a_2|$. The last
fact can possibly be explained by the bias induced by the AR part in
the FEXP estimator of $d_i$.\\

\noindent \textit{Dependent samples.}\/  Table 4 reports the performance of the test $\tilde T_n > t_{5\%}(\hat d)$ on {\it dependent} samples as in Example 3.3,
with $a_{11} = a_{22} = 1-p, \ a_{12} = a_{21} = p $, where $p \in [0,1/2)$ is a parameter. In other words, $X_i$ are
FARIMA(0,$d_i$,0) processes with mutually correlated innovations. The asymptotic correlation coefficient $\rho$
between normalized partial sums of $X_1$ and $X_2$ is proportional to $2p(1-p)/(p^2 + (1-p)^2)$ and so $\rho $ increases monotonically from 0 to 1
as $p$ increases from 0 to $1/2$. Since $p=0$ corresponds to independent samples, the results in Table 4 can be compared with those
for $a_1 = a_2=0$ in Tables 1 and 2. It appears that both tests $T_n$ and $\tilde T_n$ perform similarly and that the long-run parameter
$\rho$ is well estimated to be zero.

The purpose of Tables 5 and 6 is twofold. Firstly, we want to evaluate the performance of $\tilde T_n$ on independent
samples. Secondly, we want to assess the robustness of the adaptive formula (\ref{qopt}) for bandwidth based on AR approximation of the short
memory part with respect to other short memory specifications. To this end, we generate a FARIMA$(3,d,0)$ process with polynomial AR function
$1+0.7 x^3$ (Table 5) and a FARIMA(0,$d$,2) process with polynomial MA function $1-(1/6) x +(1/6)x^2$ (Table 6). Together with the zero hypothesis
rejection frequencies, Tables 5 and 6 also provide the (averaged) values of the adaptive estimator of the bandwidth $q$. One can immediately
recognize that
the estimated values of $q$ in Table 5 are much greater than the corresponding values on Table 6; nevertheless the size of the
$\tilde T_n$ test is respected in both tables. One can conclude from Table 6 that the adaptive formula for $q$ works rather well even
if the FAR model (on which this formula is based) is misspecified.

\newpage
 \begin{table}
\footnotesize{
\begin{tabular}{|cc|ccccc|ccccc|ccccc|}
 \hline
 \multicolumn{ 2}{|c|}{n=1024 }  &                           \multicolumn{ 5}{|c}{$a_2$=0} &                  \multicolumn{ 5}{|c}{$a_2$=0.4} &                  \multicolumn{ 5}{|c|}{$a_2$=0.8} \\
\hline
\multicolumn{ 1}{|c|}{} & $d_1\backslash d_2$& $0$& $.1$&$.2$&$.3$&$.4$& &&&&&&&&&\\
\hline
\multicolumn{ 1}{|c|}{$a_1$=0} &$0$ &   4.1 &    &    &    &    &    &    &    &    &    &    &    &    &    &    \\
\multicolumn{ 1}{|c|}{} &$.1$ &   12 &    5.0 &    &    &    &    &    &    &    &    &    &    &    &    &    \\
\multicolumn{ 1}{|c|}{} &$.2$ &    36 &    13 &    3.8 &    &    &    &    &    &    &    &    &    &    &    &    \\
\multicolumn{ 1}{|c|}{} &$.3$ &    64 &    32 &    9.3 &    3.2 &    &    &    &    &    &    &    &    &    &    &    \\
\multicolumn{ 1}{|c|}{} &$.4$ &    84 &    55 &    26 &    7.4 &    3.6 &    &    &    &    &    &    &    &    &    &    \\
\hline
\multicolumn{ 1}{|c|}{} & $d_1\backslash d_2$& $0$& $.1$&$.2$&$.3$&$.4$&$0$& $.1$&$.2$&$.3$&$.4$& &&&&\\
\hline
\multicolumn{ 1}{|l|}{$a_1$=.4} &  $0$ &   5.8 &    &    &    &    &    4.5 &    &    &    &    &    &    &    &    &    \\
\multicolumn{ 1}{|l|}{} &$.1$ & 9.1 &    4.5 &    &    &    &    13 &    3.7 &    &    &    &    &    &    &    &    \\
\multicolumn{ 1}{|l|}{} &$.2$ &    21 &    7.4 &    4.0 &    &    &    35 &    9.6 &    4.9 &    &    &    &    &    &    &    \\
\multicolumn{ 1}{|l|}{}&$.3$ &   45 &    20 &    6.6 &    3.7 &    &    59 &    33 &    12 &    3.7 &    &    &    &    &    &    \\
\multicolumn{ 1}{|l|}{} &$.4$ &      66 &    35 &    17 &    4.2 &    3.4 &    82 &    58 &    26 &    9.1 &    2.5 &    &    &    &    &    \\
\hline
\multicolumn{ 1}{|c|}{} & $d_1\backslash d_2$& $0$& $.1$&$.2$&$.3$&$.4$&$0$& $.1$&$.2$&$.3$&$.4$& $0$& $.1$&$.2$&$.3$&$.4$ \\\hline
\multicolumn{ 1}{|l|}{$a_1$=.8} &$0$ & 3.8 &    &    &    &    &    4.9 &    &    &    &    &    2.2 &    &    &    &    \\
\multicolumn{ 1}{|l|}{}&$.1$ &   3.4 &    3.6 &    &    &    &    4.0 &    4.4 &    &    &    &    9.1 &    4.5 &    &    &    \\
\multicolumn{ 1}{|l|}{} &$.2$ &   13 &    3.7 &    3.9 &    &    &    13 &    4.4 &    4.8 &    &    &    28 &    9.9 &    3.8 &    &    \\
\multicolumn{ 1}{|l|}{}&$.3$ &      26 &    10 &    3.6 &    5.2 &    &    32 &    14 &    4.0 &    4.4 &    &    56 &    31 &    11 &    3.2 &    \\
\multicolumn{ 1}{|l|}{} &$.4$ &     45 &    25 &    9.6 &    4.5 &    5.9 &    55 &    32 &    13 &    4.1 &    5.8 &    80 &    55 &    31 &    9.5 &    3.5 \\\hline
\end{tabular} }

\caption{ Frequency of rejection (in percentages) of the null hypothesis of the test $ T_n > c_{5\%}(\hat d)$. The samples are simulated following FAR(1,$d_i$) models. For fixed $a_1,a_2$, each cell contains a triangular array of dimension 5x5 corresponding to the different
parameters $(d_1,d_2)$ with $d_i \in \{0,\ .1,\ .2,\ .3,\ .4\}$ (i=1,2) and $d_1 \geq d_2$. The sample size is $n=1024$. The estimation is based on 1000 replications. }
\end{table}

\clearpage
\newpage


 \begin{table}
\footnotesize{
\begin{tabular}{|cc|ccccc|ccccc|ccccc|}
 \hline
  \multicolumn{ 2}{|c|}{n=4096 }    &                               \multicolumn{ 5}{|c}{$a_2$=0} &                  \multicolumn{ 5}{|c}{$a_2$=0.4} &                  \multicolumn{ 5}{|c|}{$a_2$=0.8} \\
\hline
\multicolumn{ 1}{|c|}{} & $d_1\backslash d_2$& $0$& $.1$&$.2$&$.3$&$.4$& &&&&&&&&&\\
\hline
\multicolumn{ 1}{|c|}{$a_1$=0} &$0$ &   5.1 &    &    &    &    &    &    &    &    &    &    &    &    &    &    \\
\multicolumn{ 1}{|c|}{} &$.1$ &   23 &    4.1 &    &    &    &    &    &    &    &    &    &    &    &    &    \\
\multicolumn{ 1}{|c|}{} &$.2$ &    59 &   20 &    4.8 &    &    &    &    &    &    &    &    &    &    &    &    \\
\multicolumn{ 1}{|c|}{} &$.3$ &    86 &    49 &   14 &   4.3&    &    &    &    &    &    &    &    &    &    &    \\
\multicolumn{ 1}{|c|}{} &$.4$ &    96 &    77 &    42 &   11 &   2.9 &    &    &    &    &    &    &    &    &    &    \\
\hline
\multicolumn{ 1}{|c|}{} & $d_1\backslash d_2$& $0$& $.1$&$.2$&$.3$&$.4$&$0$& $.1$&$.2$&$.3$&$.4$& &&&&\\
\hline
\multicolumn{ 1}{|l|}{$a_1$=.4} &  $0$ &   5.2 &    &    &    &    &    3.5 &    &    &    &    &    &    &    &    &    \\
\multicolumn{ 1}{|l|}{} &$.1$ & 11 &    4.2 &    &    &    &    20 &    5.7 &    &    &    &    &    &    &    &    \\
\multicolumn{ 1}{|l|}{} &$.2$ &    35 &    11 &    4.8 &    &    &    57 &    19 &    4.2 &    &    &    &    &    &    &    \\
\multicolumn{ 1}{|l|}{}&$.3$ &  70 &    34 &    11 &   5.2 &    &    84 &    54 &    16 &    4.0 &    &    &    &    &    &    \\
\multicolumn{ 1}{|l|}{} &$.4$ &      88 &    64 &    30 &   10&    3.3 &    95 &    78 &    48 &    14 &    3.0 &    &    &    &    &    \\
\hline
\multicolumn{ 1}{|c|}{} & $d_1\backslash d_2$& $0$& $.1$&$.2$&$.3$&$.4$&$0$& $.1$&$.2$&$.3$&$.4$& $0$& $.1$&$.2$&$.3$&$.4$ \\
\hline
\multicolumn{ 1}{|l|}{$a_1$=.8} &$0$ & 5.4 &    &    &    &    &    5.2 &    &    &    &    &    4.9 &    &    &    &    \\
\multicolumn{ 1}{|l|}{}&$.1$ &  7.4 &    3.8 &    &    &    &    9.0 &    4.3 &    &    &    &   20&    5.5 &    &    &    \\
\multicolumn{ 1}{|l|}{} &$.2$ &   25 &   9.0 &   4.3 &    &    &    30 &    10 &    4.8 &    &    &    42 &    15 &    3.8 &    &    \\
\multicolumn{ 1}{|l|}{}&$.3$ &      53 &    24 &   7.5 &    4.5 &    &   59 &    26 &   10 &   3.8 &    &   84 &   47 &   16 &    4.5 &    \\
\multicolumn{ 1}{|l|}{} &$.4$ &     75 &    44 &   21 &   6.2 &    4.2 &  82 &    50 &    25 &    8.2 &    5.3 &   95 &    79 &    50 &    16 &    4.2 \\\hline
\end{tabular} }

\caption{ The same results as Table 1 for sample size $n=4096$.}
 \end{table}

 \begin{table}
\footnotesize{
\begin{tabular}{|cc|ccccc|ccccc|ccccc|}
 \hline
  \multicolumn{ 2}{|c|}{n=4096 }    &                               \multicolumn{ 5}{|c}{$a_2$=0} &                  \multicolumn{ 5}{|c}{$a_2$=0.4} &                  \multicolumn{ 5}{|c|}{$a_2$=0.8} \\
\hline
\multicolumn{ 1}{|c|}{} & $d_1\backslash d_2$& $0$& $.1$&$.2$&$.3$&$.4$& &&&&&&&&&\\
\hline
\multicolumn{ 1}{|c|}{$a_1$=0} &$0$ &    3.2 &    &    &    &    &    &    &    &    &    &    &    &    &    &    \\
\multicolumn{ 1}{|c|}{} &$.1$ &  2.7 &    2.1 &    &    &    &    &    &    &    &    &    &    &    &    &    \\
\multicolumn{ 1}{|c|}{} &$.2$ &    2.3 &    2.0 &    1.7 &    &    &    &    &    &    &    &    &    &    &    &    \\
\multicolumn{ 1}{|c|}{} &$.3$ &   1.9 &    1.8 &    1.4 &    1.0 &    &    &    &    &    &    &    &    &    &    &    \\
\multicolumn{ 1}{|c|}{} &$.4$ &    1.8 &    1.5 &   1.0 &    0.5 &    0.3 &    &    &    &    &    &    &    &    &    &    \\
\hline
\multicolumn{ 1}{|c|}{} & $d_1\backslash d_2$& $0$& $.1$&$.2$&$.3$&$.4$&$0$& $.1$&$.2$&$.3$&$.4$& &&&&\\
\hline
\multicolumn{ 1}{|c|}{$a_1$=0.4} &$0$ &      11.2 &    &    &    &    &    5.4 &    &    &    &    &    &    &    &    &    \\
\multicolumn{ 1}{|c|}{} &$.1$ &   9.0 &    7.5 &    &    &    &    4.4 &    3.7 &    &    &    &    &    &    &    &    \\
\multicolumn{ 1}{|c|}{} &$.2$ &   7.5 &    6.3 &    5.3 &    &    &    3.6 &    3.0 &    2.7 &    &    &    &    &    &    &    \\
\multicolumn{ 1}{|c|}{} &$.3$ &   6.2 &    5.3 &    4.3 &    2.9 &    &    3.1 &    2.6 &    2.0 &    1.4 &    &    &    &    &    &    \\
\multicolumn{ 1}{|c|}{} &$.4$ &   5.3 &    4.4 &    2.9 &    1.8 &    1.0 &    2.6 &    2.0 &    1.4 &    0.8 &    0.4 &    &    &    &    &    \\
\hline
\multicolumn{ 1}{|c|}{} & $d_1\backslash d_2$& $0$& $.1$&$.2$&$.3$&$.4$&$0$& $.1$&$.2$&$.3$&$.4$& $0$& $.1$&$.2$&$.3$&$.4$ \\
\hline
\multicolumn{ 1}{|c|}{$a_1$=0.8} &$0$ &    24.8 &    &    &    &    &    22.0 &    &    &    &    &    10.2 &    &    &    &    \\
\multicolumn{ 1}{|c|}{}&$.1$ &   19.9 &    16.2 &    &    &    &    17.5 &    14.1 &    &    &    &    8.2 &    6.7 &    &    &    \\
\multicolumn{ 1}{|c|}{} &$.2$ &  16.2 &    13.4 &    11.3 &    &    &    14.2 &    11.6 &    9.6 &    &    &    6.7 &    5.7 &    4.7 &    &    \\
\multicolumn{ 1}{|c|}{}&$.3$ &  13.5 &    11.3 &    8.8 &    5.9 &    &    11.6 &    9.7 &    7.6 &    4.9 &    &    5.6 &    4.6 &    3.4 &    2.3 &    \\
\multicolumn{ 1}{|c|}{} &$.4$ &     11.3 &    8.8 &    5.8 &    3.6 &    2.1 &    9.7 &    7.5 &    5.0 &    3.0 &    1.6 &    4.5 &    3.4 &    2.2 &    1.3 &    0.6 \\
\hline
\end{tabular}
}
\caption{ The mean values on 1000 replications of $\hat q$ according to (\ref{qopt})
 for the simulations of Table 2.}
\end{table}

\begin{table}

\begin{tabular}{|rr|ccccc|ccccc|}
\hline

      &         &          \multicolumn{ 5}{|c}{$ n=1028$} &                  \multicolumn{ 5}{|c|}{$ n=4096$} \\
\hline
\multicolumn{ 1}{|c|}{} & $d_1\backslash d_2$& $0$& $.1$&$.2$&$.3$&$.4$&$0$& $.1$&$.2$&$.3$&$.4$\\
\hline
\multicolumn{ 1}{|l|}{$p=0$} &  0&   4.3 &    &    &    &    &    5.6 &    &    &    &    \\

\multicolumn{ 1}{|l|}{} &  .1&  14 &    5.4 &    &    &    &    22 &   6.4 &    &    &    \\

\multicolumn{ 1}{|l|}{} &  .2&   38 &    13 &    3.7 &    &    &    62 &    17 &    5.7 &    &    \\

\multicolumn{ 1}{|l|}{} & .3&   66 &    33 &    8.5 &    3.9 &    &    87 &    55 &    17 &    5.9 &    \\

\multicolumn{ 1}{|l|}{} &  .4&   83 &    57 &    27 &    7.5 &    3.4 &   97 &    83 &    45 &    13 &    3.8 \\




\hline
\multicolumn{ 1}{|c|}{} & $d_1\backslash d_2$& $0$& $.1$&$.2$&$.3$&$.4$&$0$& $.1$&$.2$&$.3$&$.4$\\
\hline
\multicolumn{ 1}{|l|}{$p=.15$} & 0&    5.8 &    &    &    &    &    5.8 &    &    &    &    \\

\multicolumn{ 1}{|l|}{} &  .1&  13 &    4.9 &    &    &    &    22 &    5.9 &    &    &    \\

\multicolumn{ 1}{|l|}{} &  .2&   40 &    10 &    6.2 &    &    &    64 &   21 &    6.3 &    &    \\

\multicolumn{ 1}{|l|}{} &  .3&   69 &    35 &    9.6 &    3.7 &    &    90 &    58 &    17 &    5.9 &    \\

\multicolumn{ 1}{|l|}{} &  .4&   84 &    61 &    26 &    6.0 &   3.0 &    98 &    83 &    50 &    14 &    4.1 \\




\hline
\multicolumn{ 1}{|c|}{} & $d_1\backslash d_2$& $0$& $.1$&$.2$&$.3$&$.4$&$0$& $.1$&$.2$&$.3$&$.4$\\
\hline
\multicolumn{ 1}{|l|}{$p=.35$} & 0&   5.4 &    &    &    &    &    4.9 &    &    &    &    \\

\multicolumn{ 1}{|l|}{} &  .1&   17 &    4.5 &    &    &    &    30 &    5.3 &    &    &    \\

\multicolumn{ 1}{|l|}{} &  .2&   54 &    14 &    4.7 &    &    &    84 &    26 &    5.3 &    &    \\

\multicolumn{ 1}{|l|}{} &  .3&   81 &    43 &    9.0 &    3.5 &    &    98 &    76 &    20 &    4.2 &    \\

\multicolumn{ 1}{|l|}{} &  .4&  95 &    74 &    30 &    7.2 &    2.9 &    100 &   95 &   60 &    13 &    3.1 \\
\hline
\multicolumn{ 1}{|c|}{} & $d_1\backslash d_2$& $0$& $.1$&$.2$&$.3$&$.4$&$0$& $.1$&$.2$&$.3$&$.4$\\
\hline
\multicolumn{ 1}{|l|}{$p=.45$} & 0&   6.2 &    &    &    &    &    5.5 &    &    &    &    \\

\multicolumn{ 1}{|l|}{} &  .1&  43 &    3.0 &    &    &    &    84 &    4.6 &    &    &    \\

\multicolumn{ 1}{|l|}{} &  .2&  84 &    28 &    3.6 &    &    &    98 &   63 &    5.9 &    &    \\

\multicolumn{ 1}{|l|}{} &  .3&  95 &    74 &    13 &    3.9 &    &    100 &    96 &    33 &    3.1 &    \\

\multicolumn{ 1}{|l|}{} &  .4&   97 &    90 &    46 &    6.8 &    3.2 &   100 &    99 &    85 &    14 &    4.4 \\
\hline
\end{tabular} \centering

\caption{Frequency of rejection (in percentages) of the null hypothesis of the test $\tilde T_n > c_{5\%}(\hat d)$.
 The samples are simulated following the model in (\ref{modeldep}).
 For fixed $p$, each cell contains a triangular array of dimension $5\times 5$ corresponding to the different
parameters $(d_1,d_2)$ with
 $d_i \in \{0,.1,.2,.3,.4\}$ (i=1,2) and $d_2 \leq d_1$.
}

\end{table}

\begin{table}[h]
 \centering
  \begin{tabular}[h]{|c|ccccc|}
 \hline
 $d_1 \backslash d_2$ & $0$& $.1$& $.2$& $.3$& $.4$\\
\hline
 $0$&   6.1 &    &    &   &   \\
 $.1$&  7.4  & 5.6    &&& \\
 $.2$&  21 & 8.2 & 4.0     &&      \\
 $.3$&  46  & 20  & 7.5 & 6.0  &    \\
 $.4$ &  64  & 39  & 18 &  6.0 & 4.5 \\
\hline
 \end{tabular}
 \begin{tabular}[h]{|c|ccccc|}
\hline
$d_1 \backslash d_2$ & $0$& $.1$& $.2$& $.3$& $.4$\\
\hline
 $0$& 34.7&     &  & & \\
 $.1$&  28.4& 24.0&   &&    \\
$.2$&  23.9& 20.8& 18.2&    &    \\
$.3$ &  20.8& 18.3& 16.0& 13.0&    \\
$.4$ &  18.2& 15.9& 13.0&  9.0&  5.8\\
\hline
 \end{tabular}
\caption{[Left] Frequency of rejection (in percentages) of the null hypothesis of the test $\tilde T_n > c_{5\%}(\hat d)$. The two samples $X_1$ and $X_2$ are independent. $X_1$ is simulated from FARIMA(3,$d_1$,0) model with polynomial AR function $1+ .7x^3$
  and $X_2$ from FARIMA(0,$d_2$,0) model.  [Right] Adaptive estimation of the bandwidth parameter $q$. The sample size is 4096 and the statistics are evaluated from 1000 independent replications.}
\end{table}

\begin{table}[h]
 \centering
 \begin{tabular}[h]{|c|ccccc|}
 \hline
 $d_1 \backslash d_2$ & $0$& $.1$& $.2$& $.3$& $.4$\\
\hline
 $0$&  4.9 &     &     &    &    \\
 $.1$&  17&   6.7&     &    &    \\
 $.2$&  48&  16&   4.9&    &    \\
 $.3$&  72&  44&   15&  4.3&    \\
 $.4$&  88&   71&  35&  11&  4.3\\
\hline
 \end{tabular}
 \begin{tabular}[h]{|c|ccccc|}
\hline
$d_1 \backslash d_2$ & $0$& $.1$& $.2$& $.3$& $.4$\\
\hline
 $0$&   7.0   &     &    &   & \\
 $.1$&     6.1&  5.2&   &    &    \\
 $.2$&  5.3&  4.6&  4.3&   &    \\
 $.3$&  4.8&  4.2&  3.8&  2.8&    \\
 $.4$& 4.3&  3.5&  2.7&  1.8&  1.1  \\\hline
 \end{tabular}
 \caption{[Left] Frequency of rejection (in percentages) of the null hypothesis of the test $\tilde T_n > c_{5\%}(\hat d)$.  The two samples $X_1$ and $X_2$ are independent. $X_1$ is simulated from FARIMA(0,$d_1$,2) model with polynomial MA function $1-(1/6)x +(1/6)x^2$
  and $X_2$ from FARIMA(0,$d_2$,0) model.  [Right] Adaptive estimation of the bandwidth parameter $q$. The sample size is 4096 and the statistics are evaluated from 1000 independent replications.}
\end{table}

\vskip3cm
\clearpage

\section{Concluding remarks}\label{conclusion}

The paper constructs a two-sample test for comparison of long memory parameters $d_i \in [0, 1/2), i =1,2 $ of covariance-stationary
time series $X_i, i=1,2 $ with discrete time. The test statistic, $T_n$, is defined as the sum of the ratio and the reciprocal ratio of
the rescaled variance (V/S) statistics, computed for each sample,  whose asymptotic and finite-sample behavior was
studied in Giraitis et al.~(2003, 2006). Under some assumptions which involve the existence of long-run covariances and the joint convergence of
partial sums of $X_1$ and $X_2$ to a bivariate fractional Brownian motion, we derive
the asymptotic null distribution of  $T_n$. A modification $\tilde T_n$ of the test statistic $T_n$ is introduced and shown to
be asymptotically free of the long-run correlation coefficient between
the two samples. The case when $(X_1,X_2)$ form a bivariate linear process is discussed in detail.
Simulation results using FARIMA samples  with various fractional
and autoregressive parameters show that the proposed tests have a good size for most values of fractional and autoregressive/moving average parameters.
The robustness property of the test is largely  due to our choice of bandwidth according to the adaptive formula in (\ref{qopt}) which takes into account
the estimated difference of short memory spectrum of the sampled processes.
The derivation of the last formula
uses the asymptotic expansion of the HAC estimator in Abadir et al.~(2009).

\section{Appendix. Proofs and auxiliary results}\label{appendix}

Subsection 6.1 provides alternative definition of bi-fBm by explicit cross-covariance function. Subsection 6.2 contains proofs of Propositions
\ref{propS}, \ref{Test}, \ref{TTest}, \ref{power} and \ref{prop-bilinear}. Subsection 6.3 is given to the
derivation of the adaptive bandwidth formula in (\ref{qopt}).

\subsection{Covariance function of bivariate fractional Brownian motion}
\label{sec:covar-funct-bivar}

From (\ref{Bi}) and Samorodnitsky and Taqqu~(2006) we have for any $s, t \in \R$
$$
\E B_i(s) B_i(t)= \frac{1}{2}\big\{|s|^{2d_i+1} + |t|^{2d_i+1} - |t-s|^{2d_i+1}\big\}, \qquad i=1,2.
$$
The analytic expression of cross-covariance $\E B_1(s) B_2(t)$ is derived in Lavancier et al.~(2009).
It takes a different form in the cases $d_1+d_2 \ne 0 $ and $d_1+d_2 = 0$.
Let
$$
\psi(d_1,d_2) = \frac{{\rm B}(d_1+1,d_2+1)\sqrt{\cos(d_1\pi) \cos(d_2\pi)}}{\sqrt{{\rm B}(d_1+1,d_1+1) {\rm B}(d_2+1,d_2+1)}}.
$$
Let $d_1+d_2 \ne 0 \ (d_1,d_2 \in (-1/2,1/2)$. Then
\begin{eqnarray} \label{covB12}
\hskip-.6cm \E B_1(s) B_2(t)&=&\frac{\rho_W}{2}\Big\{g_{12}(s) |s|^{d_1+d_2+1} + g_{21}(t) |t|^{d_1+d_2+1} - g_{21}(t-s)|t-s|^{d_1+d_2+1}\Big\},
\end{eqnarray}
where
$$
g_{ij}(t) = \cases{g_{ij}, &$t>0$, \cr
g_{ji}, &$t<0$ \cr}
$$
and where
\begin{equation}
 g_{12}\ =\ \psi(d_1,d_2)\sin(d_1\pi)/\sin((d_1+d_2)\pi), \quad g_{21}\ =\ \psi(d_1,d_2)
\sin(d_2\pi)/\sin((d_1+d_2)\pi).\label{eq:2}
\end{equation}

\noindent In the case $d_1+d_2 = 0$,
\begin{eqnarray} \label{covB12L}
\E B_1(s) B_2(t)&=&\frac{\rho_W }{2}\Big\{g_1 \big(|s| +|t| - |t-s|\big) \\
&& \ + \ g_2 \big(t \log |t| + s \log |s|
- (t-s) \log |t-s|\big)\Big\}, \nonumber
\end{eqnarray}
where
\begin{equation}
g_1\ =\ (1/2)\psi(d_1,d_2)(\cos (\pi d_1)+ \cos(\pi d_2)), \quad
g_2\ =\ \psi(d_1,d_2)(d_2-d_1). \label{eq:3}
\end{equation}

\subsection{Proofs of the propositions}
\label{sec:proofs-propositions}

{\bf Proof of Proposition \ref{propS}.} The first relation in (\ref{Slimit}) is immediate from (\ref{Clim1}) and (\ref{Clim2}); see also
Giraitis et al.~(2006, (3.4)). Then, the second relation in (\ref{Slimit}) follows from (\ref{SS}), which
is proved below.

Assume without loss of generality that $\mu_i = \mu_j = 0$. By Assumption A$(d_1,d_2)$, there exists a constant $C$ such that for any $n, h \ge 1,$
\begin{equation}\label{variance}
\E \Big(\sum_{t=n-h+1}^n X_i(t)\Big)^2 \ \le \ C h^{1+ 2d_i}, \qquad i=1,2.
\end{equation}
In particular, $\E \bar X^2_i = O(n^{2d_i-1})$.
Let $h \ge 1$. Then
$$
\hat\gamma_{ij}(h) - \hat\gamma_{ij}^\circ (h) \ = \ - \bar X_i \bar X_j + \bar X_i \, \frac{1}{n} \sum_{t=n-h+1}^n X_j(t)
+ \bar X_j \,\frac{1}{n} \sum_{t=n-h+1}^n X_i(t) - \frac{h}{n}\bar X_i \bar X_j.
$$
Clearly, (\ref{SS}) follows from $q/n \to 0$ and
\begin{eqnarray} \label{mean}
\frac{1}{qn} |\bar X_i| \sum_{h=1}^q \Big|\sum_{t=n-h+1}^n X_j(t)\Big|&=&o_p(n^{d_i+d_j -3/2}).
\end{eqnarray}
 By (\ref{variance}) and Cauchy-Schwarz inequality,
 \begin{eqnarray*}
\E |\bar X_i| \sum_{h=1}^q \Big|\sum_{t=n-h+1}^n X_j(t)\Big|&\le& \E^{1/2} \bar X_i^2
\E^{1/2} \sum_{h=1}^q \ \Big|\sum_{t=n-h+1}^n X_j(t)\Big|^2 \\
&\le&C n^{d_i- 1/2} \left(\sum_{h=1}^q h^{1 + 2d_j}\right)^{1/2}\\
&=& qn\  o(n^{d_i +d_j -3/2} )
 \end{eqnarray*}
since $q= o(n)$. This proves (\ref{mean}). \hfill $\Box $

\bigskip

\noindent {\bf Proof of Proposition \ref{Test}.} Both parts (i) and (ii) follow from the joint convergence
\begin{eqnarray} \label{UU}
\left((q/n)^{2d_1} V_1/S_{11,q}, (q/n)^{2d_2} V_2/S_{22,q}\right)&\rightarrow_{\rm law}&(U_1,U_2),
\end{eqnarray}
with $U_1,U_2$ as in (\ref{Uidef}). The last relation follows similarly as in Giraitis et al.~(2006, Lemmas 7.1 and 7.2).
\hfill $\Box $

\medskip

\noindent {\bf Proof of Proposition \ref{TTest}. } (i) Recall from Remark \ref{rho} that $\rho = \rho_W$. From (\ref{Slimit}), (\ref{tildeS}) and Assumptions A$(d_1,d_2)$ and B$(d_1,d_2)$ (with $d_1=d_2=d$) it follows that
\begin{eqnarray*}
q^{-2d}\tilde S_{11,q}&=&q^{-2d}S_{11,q} - \frac{(q^{-2d}S_{12,q})^2}{q^{-2d}S_{22,q}} \\
&\rightarrow_p&c_{11} - \frac{c^2_{12}}{c_{22}} \ = \ c_{11}(1-\rho^2)
\end{eqnarray*}
and
\begin{eqnarray*}
&&n^{-d-(1/2)}\left(\sum\nolimits_{t=1}^{[n\tau]} (\tilde X_1(t)- \bar{\tilde X}_1),
\sum\nolimits_{t=1}^{[n\tau]} (X_2(t)- \bar X_2) \right) \\
&&\rightarrow_{\rm fdd} \
\left(\sqrt{c_{11}} ( B_{1}(\tau) - \rho B_{2}(\tau)), \sqrt{c_{22}} B_{2}(\tau)\right) \\
&&=_{\rm fdd} \
\left(\sqrt{c_{11}(1-\rho^2)} \hat B_{1}(\tau), \sqrt{c_{22}} B_{2}(\tau)\right),
\end{eqnarray*}
where
$(\hat B_{1}(\tau) = ( B_{1}(\tau) - \rho B_{2}(\tau))/\sqrt{1-\rho^2})$ is a fBm
independent of $(B_{2}(\tau))$; see Remark \ref{rho} . Let $(\hat B^0_1(\tau) = \hat B_1(\tau) - \tau \hat B_1(1), \tau \in [0,1]), \
(B^0_2(\tau) = B_2(\tau) - \tau B_2(1), \tau \in [0,1]) $ be respective fractional Brownian bridges.

These relations together with (\ref{Slimit}) imply similarly as in Giraitis et al.~(2006) that
\begin{eqnarray}
n^{-2d}\tilde V_1&\rightarrow_{\rm law}&
c_{11}(1-\rho^2)\left(\int_0^1 (\hat B^0_{1}(\tau))^2 {\d}\tau
- \Big(\int_0^1 \hat B^0_{1}(\tau) {\d}\tau\Big)^2\right), \label{C1} \\
n^{-2d}V_2&\rightarrow_{\rm law}&c_{22}\left(
\int_0^1 (B^0_{2}(\tau))^2 {\d}\tau
- \Big(\int_0^1 B^0_{2}(\tau) {\d}\tau\Big)^2\right), \label{C2}\\
q^{-2d}\tilde S_{11,q}&\rightarrow_{\rm law}&
c_{11}(1-\rho^2), \label{C3} \\
q^{-2d}S_{22,q}&\rightarrow_{\rm law}&c_{22}\label{C4}
\end{eqnarray}
as $n, q, n/q \to \infty$, as well as the {\it joint} convergence of the four quantities
in (\ref{C1})-(\ref{C4}). Since the limits in (\ref{C1})-(\ref{C4}) are a.s. strictly positive and $(\hat B_{1}(\tau))$ is independent
of $(B_2(\tau))$, this proves (\ref{Tlimit}) and part (i).

\smallskip

\noindent (ii) From (\ref{Slimit}), (\ref{tildeS}) we have
\begin{eqnarray} \label{qtildeS}
q^{-2d_1} \tilde S_{11,q}&=&q^{-2d_1} S_{11,q} - \frac{\left(q^{-d_1-d_2}S_{12,q}\right)^2}{q^{-2d_2}S_{22,q}}\nonumber \\
&\rightarrow_p&c_{11} - \frac{c^2_{12}}{c_{22}} \ = \ c_{11}(1-\rho^2).
\end{eqnarray}
From (\ref{Slimit}) and Assumption B$(d_1,d_2)$ we obtain that
\begin{eqnarray} \label{tildeX1conv}
&&\frac{1}{n^{d_1+1/2}}\sum_{t=1}^{[n\tau]} (\tilde X_1(t)- \bar{\tilde X}_1) \nonumber \\
&&= \ \frac{1}{n^{d_1+1/2}}\sum_{t=1}^{[n\tau]} (X_1(t)- \bar{X}_1)
- \frac{q^{-d_1-d_2} S_{12,q}}{q^{-2d_2} S_{22,q}} \Big(\frac{q}{n}\Big)^{d_1-d_2} \frac{1}{n^{d_2+1/2}} \sum_{t=1}^{[n\tau]}
(X_2(t) - \bar X_2) \nonumber \\
&&\rightarrow_{\rm fdd} \
\sqrt{c_{11}} B_{1}(\tau).
\end{eqnarray}
Using similar arguments as in part (i), from (\ref{qtildeS}) and (\ref{tildeX1conv}) we get
\begin{equation} \label{d12}
\left( (q/n)^{2(d_1-d_2)} \frac{\tilde V_1/ \tilde S_{11,q}}{V_2/ S_{22,q}},
(q/n)^{2(d_2-d_1)} \frac{V_2/S_{22,q}}{\tilde V_1/ \tilde S_{11,q}}
\right)
\rightarrow_{\rm law} \left( \frac{U_1}{(1-\rho^2)U_2}, \frac{(1-\rho^2) U_2}{U_1} \right),
\end{equation}
where $U_i, i=1,2$ are defined in (\ref{Uidef}). Clearly, (\ref{d12}) implies (\ref{Tinfty}) and part (ii).

\smallskip

\noindent (iii) In this case, (\ref{qtildeS}) is again valid but (\ref{tildeX1conv}) must be changed to
\begin{eqnarray}
&&\frac{q^{d_2-d_1}}{n^{d_2+1/2}}\sum_{t=1}^{[n\tau]} (\tilde X_1(t)- \bar{\tilde X}_1) = \nonumber \\
&&= \ \ \frac{(q/n)^{d_2-d_1}}{n^{d_1+1/2}}\sum_{t=1}^{[n\tau]} (X_1(t)- \bar{X}_1)
- \frac{q^{-d_1-d_2} S_{12,q}}{q^{-2d_2} S_{22,q}}\frac{1}{n^{d_2+1/2}} \sum_{t=1}^{[n\tau]}
(X_2(t) - \bar X_2) \nonumber \\
&&\rightarrow_{\rm fdd} \
- \frac{c_{12}}{\sqrt{c_{22}}} B_{2}(\tau). \label{tildeX2conv}
\end{eqnarray}
From (\ref{tildeX2conv}) we obtain
$$
\frac{q^{2(d_2-d_1)}}{n^{2d_2}} \tilde V_1 \ \rightarrow_{\rm law} \
c_{11} \rho^2\left(\int_0^1 (B^0_2(\tau))^2 {\d}\tau
- \Big(\int_0^1 B^0_2(\tau) {\d}\tau\Big)^2\right).
$$
Combining this result with (\ref{qtildeS}) and the convergences in (\ref{C2}), (\ref{C4}), with $d= d_2$,
one obtains
$$
\frac{\tilde V_1/\tilde S_{11,q}}{V_2/ S_{22,q}} \ \rightarrow_{\rm law} \ \frac{\rho^2}{1-\rho^2},
$$
proving (\ref{Trho}). \hfill $\Box$

\vskip.3cm

\noindent {\bf Proof of Proposition \ref{power}.} We shall prove the second inequality in (\ref{A00}) only since
the first one can be proved analogously.
We shall use the following elementary inequalities: for any r.v.  $\xi, \eta \ge 0 $ and any $x >0$
\begin{eqnarray}\label{A}
&&\P \big(\frac{\xi}{\eta} \le x\big) \le \P(\xi \le \sqrt{x}) + \P \big(\eta > \frac 1{\sqrt{x}}\big), \\
&&\P (\xi \eta \le x)\le \P(\xi \le \sqrt{x}) + \P (\eta \le \sqrt{x}), \quad
\P (\xi  - \eta \le x)\le \P(\xi \le 2x) + \P (\eta > x). \nonumber
\end{eqnarray}
Denote
\begin{eqnarray*}
&&\tilde \xi_1 = n^{-2d_1} \tilde V_1, \quad \xi_1 = n^{-2d_1} V_1, \quad
\quad \xi_2 = q^{-2d_2} S_{22,q}, \quad  \tilde \xi_3 = q^{-2d_1} \tilde S_{11,q}, \\
&&\xi_3 = q^{-2d_1} S_{11,q},
 \quad \xi_4 = n^{-2d_2} V_2,
\quad x = a \left(\frac{q}{n}\right)^{2(d_1-d_2)}.
\end{eqnarray*}
Then
\begin{eqnarray} \label{A1}
\P (\tilde T_n \le a)
&\le&\P \Big( \frac{\tilde \xi_1 \xi_2 }{\tilde \xi_3 \xi_4} \le x\Big) \\
&\le&\P\big(\tilde \xi_1 \xi_2 \le x^{1/2}\big)
+ \P\big(\tilde \xi_3 \xi_4 > x^{-1/2} \big) \nonumber \\
&\le&\P\big(\tilde \xi_1 \le x^{1/4}\big)
+ \P\big(\xi_2 \le x^{1/4} \big)
+  \P\big(\tilde \xi_3 > x^{-1/4}\big)
+ \P\big(\xi_4 > x^{-1/4}\big).  \nonumber
\end{eqnarray}
Next, using the inequality $\tilde V_1 \ge (1/2) V_1 - \hat \beta^2 V_2  = (1/2) V_1 - \frac{\hat \rho^2 S_{11,q}}{S_{22,q}} V_2 $ and the facts
that $|\hat \rho| \le 1, \ V_2 \ge 0$, we get
\begin{equation} \label{tildexi}
\tilde \xi_1 \ \ge \  (1/2)\xi_1 - \eta_1,
\end{equation}
where
\begin{eqnarray*}
\eta_1&=&\frac{S_{11,q}}{S_{22,q}}(n^{-2d_1} V_2) \
= \  x\frac{\xi_3 \xi_4}{\xi_2}.
\end{eqnarray*}
Relations (\ref{tildexi}) and (\ref{A}) yield
\begin{eqnarray} \label{A2}
\P\big(\tilde \xi_1 \le x^{1/4}\big)
&\le&\P\big((1/2)\xi_1 - \eta_1
 \le x^{1/4}\big)\\
&\le&\P (\xi_1 \le 4 x^{1/4}) +  \P (\eta_1 > x^{1/4}) \nonumber  \\
&=&\P (\xi_1 \le 4 x^{1/4}) +  \P (x \frac{\xi_3 \xi_4}{\xi_2} > x^{1/4}) \nonumber \\
&\le&\P (\xi_1 \le 4 x^{1/4}) +  \P (\xi_3 > x^{-1/4}) + \P (\xi_4 > x^{-1/4}) + \P (\xi_2 \le x^{1/4}). \nonumber
\end{eqnarray}
Combining (\ref{A1}) and (\ref{A2}) and using $\tilde S_{11,q} \ge S_{11,q}$, see (\ref{tildeS}), we obtain
\begin{eqnarray} \label{A3}
\P (\tilde T_n \le a)
&\le&\P\big(\xi_1 \le 4x^{1/4}\big)
+ 2\P\big(\xi_2 \le x^{1/4} \big)
+  2\P\big(\xi_3 > x^{-1/4}\big)
+ 2\P\big(\xi_4 > x^{-1/4}\big).
\end{eqnarray}
From Assumption A$(d_1,d_2)$ and  (\ref{rest}) we obtain
\begin{eqnarray} \label{A4}
\E S_{11,q}&=&\frac{1}{q+1} \E (\sum_{t=0}^q X_1(t) )^2 + \E \sum_{|h| \le q} (1 - \frac{|h|}{q+1}) (\hat \gamma_{11} (h) -
\hat \gamma^\circ_{11}(h)) \\
&\le&K_2 q^{2d_1}, \nonumber \\
\E V_2&\le&K_3 n^{2d_2}  \nonumber
\end{eqnarray}
for some constants $K_2, K_3 $ independent of $n, q $, implying $\E \xi_3 \le K_2, \ \E \xi_4 \le K_3 $. From (\ref{A3}), (\ref{A4}), the Markov
inequality, and assumption (\ref{A0}), the statement
of the proposition easily follows.
\hfill $\Box$

\vskip.3cm

\noindent {\bf Proof of Proposition \ref{prop-bilinear}.} With exception of (\ref{Clim2}), all other facts
in the statement of the proposition follow similarly or using the argument developed in Davydov~(1970), Bru\v zait\.e and Vai\v ciulis~(2005),
Giraitis et al.~(2006) and other papers. In particular,
the joint convergence of partial sums of $(X_1, X_2)$ can be proved by using the
scheme of discrete stochastic integrals in Surgailis~(2003). See also Chung~(2002, proof of Theorem 1).

Let us prove the convergence of empirical long-run covariances in (\ref{Clim2}) or, equivalently, in (\ref{Slimit}). Denote
\begin{equation}
X_{ij}(t) \ = \ \sum_{k=0}^\infty \psi_{ij}(k) \xi_j(t-k), \quad i,j=1,2 \label{Xij}
\end{equation}
It suffices
to show the convergence of the HAC estimates of long-run covariances $c_{ij, i'\!j'}$ of components $X_{ij}$ and $X_{i'\!j'}$ in (\ref{Xij}), for any pairs
$(i,j), (i',j') \in \{1,2\} \times \{1,2 \}$; more precisely, to show that
\begin{eqnarray}
q^{-d_{ij}-d_{i'j'}} S_{ij, i'j', q}&\rightarrow_p&c_{ij,i'\!j'} \qquad (i,j, i',j'=1,2), \label{Slimit1}
\end{eqnarray}
where $S_{ij, i'j', q}$ is defined as in (\ref{Sdef}) with $\hat \gamma_{ij}(h)$ replaced by the empirical
covariance $\hat \gamma_{ij, i'j'}(h)$ between observations $X_{ij}(t), t=1, \cdots, n$ and
$X_{i'j'}(t), t=1, \cdots, n$.

Fix $i,j,i',j' \in \{1,2\}$ and denote
$X(t) = X_{ij}(t), \ X'(t) = X_{i'j'}(t), \ d = d_{ij}, \ d' = d_{i'j'}, $
$$
S_q = S_{ij, i'j',q}, \quad \gamma(h) = \gamma_{ij,i'j'}(h), \quad \hat \gamma (h) = \hat \gamma_{ij, i'j'}(h), \quad c = c_{ij, i'j'},
$$
$\psi(k) =
\psi_{ij}(k), \ \psi'(k) = \psi_{i'j'}(k), \
\xi(s) = \xi_j(s), \ \xi'(s) = \xi_{j'}(s), \ \tilde \rho_\xi = \rho_{\xi, jj'}$ for short. Write $S_{q} = S'_{q} + S''_{q}, $ where
\begin{eqnarray*}
S'_{q}&=&\sum_{h=-q}^q \left(1-\frac{|h|}{q+1}\right) \tilde \gamma(h), \qquad S''_{q}\ = \
\sum_{h=-q}^q \left(1-\frac{|h|}{q+1}\right) (\hat \gamma (h) - \tilde \gamma(h)),
\end{eqnarray*}
where
$$
\tilde \gamma(h)\ =\ n^{-1}\cases{\sum_{t=1}^{n-h}X(t) X'(t+h), &$h\ge 0$, \cr
\sum_{t=1-h}^{n} X(t) X'(t+h),&$h\le 0$\cr}.
$$
is the empirical covariance from noncentered observations; c.f. (\ref{hatgamma}).
Then (\ref{Slimit1}) follows from
\begin{eqnarray}
q^{-d-d'} S'_{q}&\rightarrow_p&c, \qquad S''_{q} \ = \ o_p(q^{d+d'}). \label{Slimit2}
\end{eqnarray}
In the subsequent proof of the first relation of (\ref{Slimit2}), we first assume $d >0, d' >0$. Split $S'_q
= \tilde \gamma(0) + \sum_{h=-q}^{-1} \left(1-\frac{|h|}{q+1}\right) \tilde \gamma(h) + \sum_{h=1}^q \left(1-\frac{|h|}{q+1}\right) \tilde \gamma(h)$. Here, the last two sums can be treated similarly and $\tilde \gamma(0) = O_p(1) $ is negligible.
Consider
\begin{eqnarray*}
\sum_{h=1}^q \left(1-\frac{h}{q+1}\right) \tilde \gamma(h)
&=&\sum_{i=1}^3 U_i,
\end{eqnarray*}
where
\begin{eqnarray}
U_1&=&
\tilde \rho_\xi \sum_{s} \sum_{h=1}^q \left(1-\frac{h}{q+1}\right) \frac{1}{n} \sum_{t=1}^{n-h} \psi(t-s)\psi'(t+h-s), \nonumber \\
U_2&=&
\sum_{s} \eta_s \sum_{h=1}^q \left(1-\frac{h}{q+1}\right) \frac{1}{n} \sum_{t=1}^{n-h} \psi(t-s)\psi'(t+h-s), \label{U2} \\
U_3&=&\sum_{s \ne s'} \sum_{h=1}^q \left(1-\frac{h}{q+1}\right)
\frac{1}{n} \sum_{t=1}^{n-h} \psi(t-s)\psi'(t+h-s')
\xi(s)\xi'(s'), \label{U3}
\end{eqnarray}
where $\eta_s = \xi(s)\xi'(s) - {\E} \xi(s) \xi'(s) = \xi(s) \xi'(s) - \tilde \rho_\xi$,
the sums $\sum_s$ and $ \sum_{s\ne s'}$ are taken over all $s\in {\Z} $ and $s,s' \in {\Z}, s \ne s'$, respectively,
and where we put $\psi(t) = \psi'(t) = 0 \ (t<0)$.

First, consider the (nonrandom) term $U_1$.
Using the asymptotics of $\psi$ and $\psi'$ and the dominated convergence theorem, we easily obtain that, as $q \to \infty,\ n \to \infty,\ q/n \to 0$,
\begin{eqnarray*}
U_1&=&\tilde \rho_\xi \sum_{h=1}^q \left(1-\frac{h}{q+1}\right) \Big(\frac{1}{n} \sum_{t=1}^{n-h} 1 \Big)
\sum_{k=0}^\infty \psi (k) \psi'(h+k) \\
&\sim&\tilde \rho_\xi \,\alpha \,\alpha'\sum_{h=1}^q \int_0^\infty k^{d-1} (h+k)^{d'-1} {\d}k \\
&\sim&\tilde \rho_\xi \, \alpha\, \alpha' {\rm B}(d, 1-d-d') \sum_{h=1}^q h^{d+ d' -1} \\
&\sim& \ \frac 12 c q^{d + d'},
\end{eqnarray*}
where
$
c = 2\tilde \rho_\xi \, \alpha \, \alpha' \, {\rm B}(d, 1-d-d')/(d+ d'). $
Then, the first relation in (\ref{Slimit2}) follows from
\begin{equation}
q^{-d-d'}U_i \ = \ o_p(1) \qquad (i=2,3). \label{Slimit3}
\end{equation}

To estimate $U_2$, we use the fact that $(\eta(s), s \in \Z)$ are i.i.d.r.v.'s, the
well-known inequality ${\E} |\sum_i M_i|^p \le 2 \sum_i {\E} |M_i|^p $ for independent zero mean
random variables $M_i$ with ${\E}|M_i|^p < \infty $
 and $1\le p \le 2$ (see Bahr and Ess\'een~(1965)), the fact ${\E} |\eta_s|^p = C_p< \infty $ for some
$p \in (1,2) $ and the Minkowski inequality. Using these facts, we obtain
\begin{eqnarray*}
{\E}|U_2|^p&\le&
2C_p\sum_{s} \bigg( \sum_{h=1}^q \frac{1}{n} \sum_{t=1}^{n} |\psi(t-s)\psi'(t+h-s)| \bigg)^p \\
&\le&C\bigg(\sum_{h=1}^q \frac{1}{n} \sum_{t=1}^n \Big( \sum_s |\psi(t-s) \psi'(t+h-s)|^p \Big)^{1/p} \bigg)^p \\
&\le&C \bigg(\sum_{h=1}^q \Big(\sum_{s=0}^\infty s_+^{p(d-1)} (h+s)^{p(d'-1)} \Big)^{1/p} \bigg)^p \\
&\le&C \bigg(\sum_{h=1}^q \Big(h^{p(d+d'-2)+1} \Big)^{1/p} \bigg)^p \\
&\le&C \Big(\sum_{h=1}^q h^{(d+d'-2)+(1/p)} \Big)^p \ \le \ C q^{p((1/p)+d+d'-1)}
\end{eqnarray*}
and therefore ${\E}^{1/p} |U_2|^p = O(q^{d+d' +(1/p)-1}) = o(q^{d+d'})$, as $p> 1$, proving
(\ref{Slimit3}) for $i=2$.

\smallskip

Next, consider $U_3$. Using the fact that $\sum_{s=1}^\infty s^{d-1} (t+s)^{d'-1} \le C t^{d+d'-1} \ (t\ge 0)$, we obtain
\begin{eqnarray}
{\E}U^2_3&\le&2\sum_{s \ne s'} \bigg( \sum_{h=1}^q \left(1-\frac{h}{q+1}\right)\frac{1}{n} \sum_{t=1}^{n-h} \psi(t-s)\psi'(t+h-s') \bigg)^2\nonumber\\
&\le&\frac{C}{n^2}\sum_{t,t'=1}^n \sum_{h,h'=1}^q \sum_{s,s'}
|\psi(t-s)\psi'(t+h-s')\psi(t'-s)\psi'(t'+h'-s')|\label{rest0}\\
&\le&\frac{C}{n^2} \sum_{t,t'=1}^n \sum_{h,h'=1}^q |t-t'|_+^{2d-1} |t-t'+ h-h'|_+^{2d'-1}\nonumber \\
&\le&\frac{C q}{n} \sum_{|t|\le n} \sum_{|h|\le q} |t|_+^{2d-1} |t + h|_+^{2d'-1} \ \le \ C(J_1 +J_2),\nonumber
\end{eqnarray}
where
\begin{eqnarray*}
J_1&=&(q/n)\sum_{|t|\le 2q}|t|_+^{2d-1} \sum_{|h|\le 3q} |h|_+^{2d'-1} \ \le \ C(q/n)q^{2(d+d')} \ = \ o(q^{2(d+d')}), \\
J_2&=&(q^2/n)\sum_{2q <|t|\le n} |t|^{2d+2d'-2}\
\le \ C(q^2/n)\cases{n^{2d+2d' -1}, &$d+d'> 1/2$\cr
q^{2d+2d'-1}, &$d+d' < 1/2$\cr
\log (n/q), &$d+d'=1/2$\cr}
\end{eqnarray*}
and so $J_2 = o(q^{2(d+d')}) $ as $q, n, n/q \to \infty$ in all three cases (in the last case $d+d'=1/2$ this follows from the fact that
$x \to 0$ \  entails $x \log (1/x) \to 0$).

This proves the first relation of (\ref{Slimit2}) for $d>0, d' >0$.
\smallskip

Consider now  the case $d=d'=0$  we want to prove that
$$
S'_{q}\ \rightarrow_p \ c,
$$
where
\begin{eqnarray}
c&=&\lim_{n\to \infty} n^{-1} {\E} \bigg(\sum_{t=1}^n X(t)\bigg) \bigg(\sum_{s=1}^n X'(s) \bigg)\nonumber  \\
&=&\tilde \rho_\xi \, \lim_{n\to \infty} n^{-1} \sum_u \sum_{t=1}^n \psi(t-u) \sum_{s=1}^n \psi'(s-u) \nonumber \\
&=&\tilde \rho_\xi \, \alpha \, \alpha'. \label{rest2}
\end{eqnarray}
By writing $S'_q = {\E} S'_q + (S'_q - {\E} S'_q) $, the convergence ${\E} S'_q \to c $ follows
similarly as in (\ref{rest2}) above. Relation $ S'_q - {\E} S'_q = o_p(1)$ can be shown similarly to (\ref{Slimit3}), i.e., by splitting
$ S'_q - {\E} S'_q $ into ``diagonal'' and ``off-diagonal'' parts in the quadratic form in noise variables.
Consider the ``diagonal'' part $U_2$ in (\ref{U2}). Then
\begin{eqnarray*}
{\E}|U_2|^p&\le&
2C_p\sum_{s} \bigg( \sum_{h=1}^q \frac{1}{n} \sum_{t=1}^{n} |\psi(t-s)\psi'(t+h-s)| \bigg)^p
\ \le \ C(W_1 + W_2),
\end{eqnarray*}
where
\begin{eqnarray*}
W_1&=&\sum_{s=1}^n \bigg( \sum_{h=1}^q \frac{1}{n} \sum_{t=1}^{n} |\psi(t-s)\psi'(t+h-s)| \bigg)^p\\
&\le&Cn^{-p} \sum_{s=1}^n \bigg(\sum_{t=1}^\infty |\psi(t-s)| \bigg)^p \ = \ O(n^{1-p}) \ = o(1)
\end{eqnarray*}
since $p>1$, while
\begin{eqnarray*}
W_2&=&\sum_{s=0}^\infty \bigg( \sum_{h=1}^q \frac{1}{n} \sum_{t=1}^{n} |\psi(t+s)\psi'(t+h+s)| \bigg)^p\\
&\le&Cn^{-p}\sum_{s=0}^\infty \bigg(\sum_{t=1}^{n} |\psi(t+s)| \bigg)^p \\
&\le&C \bigg(n^{-1} \sum_{t=1}^n \bigg(\sum_{s=0}^\infty |\psi(t+s)|^p \bigg)^{1/p} \bigg)^p \ = \ o(1)
\end{eqnarray*}
where we used the Minkowski inequality and the dominated convergence theorem to get $o(1)$, in view of the
fact that $\sum_{s=0}^\infty |\psi(t+s)|^p$ is bounded in $t$ and tends to zero as $t \to \infty $.

Consider the ``off-diagonal'' term $U_3$ in (\ref{U3}). Noting that, for
fixed $h, h'$, the sum in (\ref{rest0}) over all
$t, t', s, s' \in {\Z}$ is bounded by a constant independent of $h, h'$, we get ${\E} U^2_3 \le C(q/n)^2 = o(1)$.

This proves the first relation of (\ref{Slimit2}) for $d>0, d' >0$ and
$d =d'=0$.

Let us prove the second relation in (\ref{Slimit2}). It follows from
\begin{equation}
\sum_{|h| \le q} {\E} |\hat \gamma(h) - \tilde \gamma(h)| = o(q^{d+d'}). \label{rest}
\end{equation}
Using definitions of $\hat \gamma (h), \tilde \gamma (h)$, the Cauchy-Schwartz inequality and (\ref{Clim1}),
for $h\ge 0$ one obtains
\begin{eqnarray}
&&\left({\E} |\hat \gamma(h) - \tilde \gamma(h)|\right)^2 \nonumber \\
&&\le \ {\E}\left(\bar {X}\right)^2 {\E}\bigg(\frac{1}{n}\sum_{t=1}^{n-h} X'(t+h)\bigg)^2 +
{\E}\left(\bar{X'}\right)^2 {\E}\bigg(\frac{1}{n}\sum_{t=1}^{n-h} X(t)\bigg)^2 +
{\E}\left(\bar{X}\right)^2 {\E}\left(\bar{X'}\right)^2 \nonumber \\
&&\le \ Cn^{2d + 2d'- 2} \label{rest1}
\end{eqnarray}
and so (\ref{rest}) reduces to $ Cq n^{d+d'-1} = o(q^{d+d'}) $ which is a consequence of $d+d'< 1$ and $q/n \to 0$.

\medskip
This concludes the proof of
(\ref{Slimit2}) for $d>0, d'>0$ and
$d =d'=0$. The cases $d>0 = d'$ and $d=0< d'$ can be treated in a similar way. Proposition \ref{prop-bilinear} is proved.
\hfill $\Box$

\subsection{Derivation of the adaptive bandwidth}
\label{sec:deriv-adapt-bandw}
The aim of this section is to derive the adaptive bandwidth formula used in our simulations, viz.
\begin{equation}\label{qopt}
\hat q \ = \ 0.3 |\hat I|^{1/2}\cases{n^{1/(3+4\hat d)}, & if $\hat d<1/4$,\cr
n^{1/2-\hat d}, & if $\hat d>1/4$,}
\end{equation}
where $\hat d = (\hat d_1 + \hat d_2)/2$ is an estimator of the (common) long memory parameter $d$,
\begin{equation}\label{hatIi}
\hat I\ =\ \int_0^{\pi} \Big(\frac{\hat g_1(x)}{\hat g_1(0)} - \frac{\hat g_2(x)}{\hat g_2(0)}\Big)\,\frac{{\d}x}{x^{2\hat d} \sin^2(x/2)},
\end{equation}
and where $\hat g_i$ is an estimator of the short memory part $g_i(x) = f_i(x)/|x|^{2d_i} $ of the
spectral density $f_i$ of $X_i, \ i=1,2$. In this paper, $\hat g_i$ is
the spectral density of the best AR approximation of $g_i$ which is computed following the two step procedure in Ray
and Crato~(1996). Namely, we first estimate $d_i$
and then we fit an AR process to $(1-L)^{\hat d_i} X_i$ using
the BIC criterion.

From Abadir et al.~(2009, Theorem 2.1) under similar assumptions on $X_i$ as in Section~\ref{test-construction}
we have the following expansion
of $S_{ii,q}$: for $0< d_i < 1/4$,
\begin{eqnarray} \label{distaso}
q^{-2d_i}S_{ii,q}&=&c_{ii} + (q/n)^{1/2}g_i(0)(Z_{ni} + o_p(1)) + q^{-1-2d_i}g_i(0)(B_i + o_p(1)),
\end{eqnarray}
where $Z_{ni} \rightarrow_{\rm law} Z_i \sim N(0, v(d_i))$,
\begin{eqnarray*}
v(d_i)&=&8 \pi \int_0^\infty \Big(\frac{\sin(x/2)}{x/2}\Big)^4 x^{-4d_i} {\d}x, \\
B_i&=&\int_0^\infty \Big( \frac{g_i(x)}{g_i(0)\sin^2(x/2)}{\bf 1}_{\{0<x<\pi\}} - \frac{1}{(x/2)^2} \Big) \frac{{\d}x}{x^{2d_i}}
\end{eqnarray*}
and where $g_i(x) = f_i(x) |x|^{2d_i} $ is the short memory component of the spectral density $f_i$ of $X_i$, which is assumed
to be continuous at $x=0$ and $g_i(x) = g_i(0) + O(x^2), \ x \to 0, \ g_i(0)>0$. Note that the long-run variance $c_{ii}$ is related to $g_i(0) $ by
\begin{equation} \label{p(d)}
c_{ii} = g_i(0) p(d_i),
\end{equation}
where $p(d) = 2\Gamma (1-2d)\sin (\pi d)/d(1+2d)$ depends only on $d$.

From the form of statistic $T_n$ it is clear that $q$ must be chosen so that the ratio $c_{11}/c_{22}$ is well estimated
by $S_{11,q}/S_{22,q}$.
From (\ref{distaso}), assuming $d_1=d_2 =d $ as under the null hypothesis, we obtain
\begin{eqnarray} \label{SSratio1}
\hskip-.5cm \frac{S_{11,q}/c_{11}}{S_{22,q}/c_{22}} -1&=&(q/n)^{1/2}\frac{Z_{n1} - Z_{n2}}{p(d)} (1+ o_p(1)) + q^{-1-2d} \frac{B_1 - B_2}{p(d)}(1+ o_p(1)).
\end{eqnarray}
Therefore as $n, q, q/n \to \infty $,
\begin{eqnarray} \label{distaso1}
E \Big(\frac{S_{11,q}/c_{11}}{S_{22,q}/c_{22}} -1\Big)^2&\sim&\frac{1}{p(d)^2}\left((q/n) \E \big(Z_{1} - Z_{2}\big)^2
+ q^{-2(1+2d)} I^2\right),
\end{eqnarray}
since for $d_1=d_2 = d$, we have
$B_1 - B_2 = I$, where
$$
I \ = \ \int_0^\pi \Big(\frac{g_1(x)}{g_1(0)} - \frac{g_2(x)}{g_2(0)}\Big) \frac{{\d}x}{x^{2d} \sin^2(x/2)},
$$
c.f. (\ref{hatIi}). Minimizing the right-hand side of (\ref{distaso1}) with respect to $q$, we obtain
\begin{equation} \label{q1}
q\ = \ K_1(d)|I|^{2/(3+4d)} n^{1/(3+4d)},
\end{equation}
where $K_1(d)$ depends on $d$ and $\E \big(Z_{1} - Z_{2}\big)^2 $.
Numerical computation of the function $K_1(d)$ reveals that it is well approximated by the constant value 0.3
on the interval $(0, 1/4)$ except for the case when $d$ is close to $1/4$ and then $K_1(d)$ diverges but then also the approximations
in (\ref{distaso}) and (\ref{distaso1}) are less accurate. Therefore we choose to replace $K_1(d)$ in (\ref{q1}) by 0.3
on the whole interval $d\in (0,1/4)$ in order that the test is not too conservative. For similar reasons, we replace
the exponent of $|I|$ in (\ref{q1}) by 1/2, since otherwise the test turns out to be too conservative
for small values of $d$ when $(X_1)$ and $(X_2)$ have very different short memory parts (or high values of $|I|$). The result
of these replacements is $q = 0.3|I|^{1/2} n^{1/(3+4d)}$, c.f. (\ref{qopt}).

Next, let us turn to the case $1/4 < d< 1/2$. From Abadir et al.~(2009, Theorem 2.1) we obtain that for $1/4 < d_i < 1/2$
\begin{eqnarray} \label{distaso2}
q^{-2d_i}S^\circ_{ii,q}&=&c_{ii} + (q/n)^{1-2d_i} g_i(0)(\tilde Z_{ni} + o_p(1)) + q^{-1-2d_i}g_i(0)(B_i + o_p(1)),
\end{eqnarray}
where $\tilde Z_{ni} \rightarrow_{\rm law} \tilde Z_i $ and $\tilde Z_i$ is a (non-Gaussian) r.v. whose distribution depends only on $d_i$.
From Proposition 2.4 (\ref{SS}) we have that
\begin{eqnarray} \label{SS1}
q^{-2d_i}(S_{ii,q} - S^\circ_{ii,q})&=&-2(q/n)^{1-2d_i} g_{i}(0) (Y^2_{ni} + o_p(1)),
\end{eqnarray}
where
\begin{equation} \label{Yn}
Y_{ni} = g_i(0)^{-1/2} n^{1/2 - d_i} \bar X_i \ \rightarrow_{\rm law} \ Y_i \sim N\Big(0, \frac{c_{ii}}{g_i(0)}\Big) = N(0, p(d_i)),
\end{equation}
see (\ref{p(d)}). Combining (\ref{distaso2})-(\ref{Yn}) and using the facts that $\E \tilde Z_i = 0, \ i=1,2 $ and $ \E Y_1^2 = \E Y^2_2 $,
 similarly as in (\ref{SSratio1}) and (\ref{distaso1})
we obtain
\begin{eqnarray*} \label{SSratio2}
\frac{S_{11,q}/c_{11}}{S_{22,q}/c_{22}} -1&=&(q/n)^{1-2d}\frac{\tilde Z_{n1} - 2Y_{n1}^2 - \tilde Z_{n2} + 2Y_{n2}^2}{p(d)}\big)(1+ o_p(1))\\
&+&q^{-1-2d} \frac{B_1 - B_2}{p(d)}(1+ o_p(1))
\end{eqnarray*}
and
\begin{eqnarray} \label{distaso3}
E \Big(\frac{S_{11,q}/c_{11}}{S_{22,q}/c_{22}} -1\Big)^2&\sim&\frac{1}{p(d)^2}\left( (q/n)^{2(1-2d)}\E J^2(d)
+ q^{-2(1+2d)} I^2\right),
\end{eqnarray}
where $J(d) = \tilde Z_{1} - 2Y^2_1- \tilde Z_{2} + 2Y_2^2 $ has a distribution depending on $d$ alone. Minimization of the
right-hand side of (\ref{distaso3}) with respect to $q$ leads to
\begin{equation} \label{q2}
q \ = \ K_2(d)|I|^{1/2} n^{1/2 - d},
\end{equation}
where $K_2(d)$ is a function of $d$. In this case, we also choose $K_2(d) = 0.3 $ for similar reasons
as in the case $d<1/4$ above. This completes our derivation of the bandwidth formula (\ref{qopt}).

\vskip2cm

\section*{Acknowledgements} The authors are grateful to Referee and Associated Editor for
comments and suggestions which helped to improve the original version of the paper.
The joint work was partially supported by the research project MATPYL
of the F\'ed\'eration de Math\'ematiques des Pays de Loire.
The research of the third author (D.S.) was partially supported by the Lithuanian State Science and Studies
Foundation grant no. T-70/09.

\section*{References}

\smallskip
\bigskip

\begin{description}
\itemsep -.04cm

\item Abadir, K., Distaso, W., Giraitis, L., 2009. Two estimators of the long-run
variance: beyond short memory. Journal of Econometrics 150, 56--70.

\item Bahr, von B., Ess\'een, C.-G., 1965. Inequalities for the $r$th absolute moment of a sum of random variables, $1\le r\le 2$.
Annals of Mathematical Statistics 36, 299--303.

\item Beran, J., Terrin, N., 1996. Testing for a change of the long-memory parameter. Biometrica
83, 627--638.

\item Bru\v zait\.e, K., Vai\v ciulis, M., 2005. Asymptotic independence of distant partial sums of linear process.
Lithuanian Mathematical Journal 45, 387--404.

\item Casas, I., Gao, J., 2008. Econometric estimation in long-range dependent volatility models: Theory and practice.
Journal of Econometrics 147, 72--83.

\item Chung, C.-F., 2002. Sample means, sample autocovariances, and linear regression of stationary multivariate long memory processes. Econometric Theory 18, 51--78.

\item Crato, N., Ray, B.K., 1996. Model selection and forecasting for long-range
dependent processes. Journal of Forecasting 15, 107--125.

\item Davydov, Yu., 1970. The invariance principle for
stationary processes. Theory Probability and Its Applications 15,
487--498.

\item Deo, R., Hurvich, C.M., Soulier, Ph., Yi Wang, 2009.
Propagation of memory parameter from durations to
counts. Econometric Theory 25 (to appear).

\item Didier, G., Pipiras,  V., 2010.  Integral representations of operator
  fractional Brownian motion. To appear in Bernoulli.

\item Giraitis, L., Kokoszka, P., Leipus, R.,  Teyssi\`ere, G.,
2003. Rescaled variance and related tests for long memory in
volatility and
  levels. Journal of Econometrics 112, 265--294.

\item Giraitis, L., Leipus, R., Philippe, A., 2006. A test for
  stationarity versus trends and unit roots for a wide class of dependent
  errors. Econometric Theory 22, 989--1029.

\item Horv\'ath, L., 2001. Change-point detection in long-memory processes. Journal of Multivariate Analysis 78, 218--234.

\item Iouditsky, A., Moulines, E., Soulier, Ph., 2001.
Adaptive estimation of the fractional differencing coefficient.
Bernoulli 7, 699--731.

\item Kwiatkowski, D., Phillips, P.C.B., Schmidt, P., Shin, Y., 1992. Testing the null hypothesis of stationarity
against the alternative of a unit root: how sure are we that
economic time series have a unit root? Journal of Econometrics
54, 159--178.

\item Lavancier, F., Philippe, A., Surgailis, D., 2009. Covariance function of vector self-similar process.
Statist. Probab. Letters 79, 2415--2421.

\item Lo, A.W., 1991. Long-term memory in stock market prices.
  Econometrica 59, 1279--1313.

\item Lobato, I., Robinson, P.M., 1998.
A nonparametric test for I(0).
Revue of Economic Studies 65, 475--495.

\item Robinson, P.M., 1994. Efficient tests of nonstationary hypotheses.
Journal of American Statistical Association 89, 1420--1437.

\item Samorodnitsky, G., Taqqu, M.S., 1994. Sta\-ble Non-Gaussian Ran\-dom Pro\-ces\-ses.
Chap\-man and Hall, New York.

\item Soofi, A.S., Wang, S., Wang, Y., 2006. Testing for long memory in the Asian foreign exchange rates.
Journal of System Science and Complexity 19, 182--190.

\item Surgailis, D., 2003. Non CLTs: U-statistics, multinomial formula and
approximations of multiple It\^o-Wiener integrals. In: Doukhan, P.,
Oppenheim, G., Taqqu, M.S. (Eds.), Theory and Applications
of Long-Range Dependence: Theory and Applications, pp. 129--142. Birkh\"auser, Boston.

\item Surgailis, D., Teyssi\`ere, G., Vai\v ciulis, M., 2008.
The increment ratio statistic. Journal of Multivariate Analysis 99, 510--541.

\item Teyssi\`ere, G., Kirman, A.P. (Eds.), 2007. Long memory in economics. Springer, Berlin - Heidelberg.

\end{description}

\end{document}